\documentclass[journal]{IEEEtran}

\usepackage{array}

\newtheorem{theorem}{Theorem} 
\newtheorem{lemma}{Lemma}

\newtheorem{definition}{Definition} 
\newtheorem{example}{Example}      
\usepackage{algorithm}
\usepackage{algorithmicx}
\usepackage{algcompatible}
\usepackage{algpseudocode}
\usepackage{graphicx}
\usepackage{textcomp}
\usepackage{color}
\usepackage{cite}
\usepackage{amsmath,amssymb,amsfonts}
\usepackage{CJKutf8}
\usepackage{CJKutf8}
\usepackage{orcidlink}
\usepackage{multirow}
\usepackage{url}

\begin{document}

\title{A Type 2 Fuzzy Set Approach for Building Linear Linguistic Regression Analysis under Multi Uncertainty}

\author{Junzo Watada $^{*1}$\orcidlink{0000-0002-3322-2086}, 
Pei-Chun Lin $^{*2}$\orcidlink{0000-0003-0735-2693},
Bo Wang $^{3}$\orcidlink{0000-0002-9264-9741}, 
Jeng-Shyang Pan $^{4}$\orcidlink{0000-0002-3128-9025}, and 
Jos\'{e} Guadalupe Flores Mu\~{n}iz $^{5}$\orcidlink{0000-0001-5318-5860}

\thanks{Junzo Watada $^{*1}$\orcidlink{0000-0002-3322-2086} is a Specially Appointed Professor in Faculty of Data Science, Shimonoseki City University, Japan; watada-ju@shimonoseki-cu.ac.jp; 
Pei-Chun Lin $^{*2}$\orcidlink{0000-0003-0735-2693} is an Associate Professor in Department of Information Engineering and Computer Science, Feng Chia University; peiclin@fcu.edu.tw;
Bo Wang $^{3}$\orcidlink{0000-0002-9264-9741} is an Associate Professor in School of Management and Engineering, Nanjing University, Nanjing, China; bowangsme@nju.edu.cn;
Jeng-Shyang Pan $^{4}$\orcidlink{0000-0002-3128-9025} is a Proffesor in Nanjing University of Information Science \& Technology, Nanjing, China, and College of Computer Science and Engineering, Shandong University of Science and Technology, Qingdao, China; jengshyangpan@gmail.com;
Jos\'{e} Guadalupe Flores Mu\~{n}iz $^{5}$\orcidlink{0000-0001-5318-5860} is a Professor and Researcher in the Centro de Investigaci\'on de Ciencias F\'isico Matem\'aticas (CICFIM) of the Universidad Aut\'onoma de Nuevo Le\'on (UANL) in Mexico}
}

\markboth{Journal of \LaTeX\ Class Files,~Vol.~14, No.~8, August~2015}%
{Shell \MakeLowercase{\textit{et al.}}: Bare Demo of IEEEtran.cls for IEEE Journals}

\maketitle

\begin{abstract}
In this paper, we propose a novel heuristic algorithm for constructing a Type-2 Fuzzy Set of the Linear Linguistic Regression (T2F-LLR) model, designed to address uncertainty and vagueness in real-world decision-making. 
We consider a practical scenario involving a cosmetic company's promotional planning across four product categories: Basic Face Care, Face Cleaning, Cosmetics, and Body Care, aimed at both male and female consumers. 
Data are collected using fuzzy linguistic questionnaires from customers and expert managers, with responses expressed using qualitative terms such as 'always', 'frequently', 'Often', 'Sometimes', and 'frequently'. These linguistic evaluations are modeled as Type-2 Fuzzy Set of Linear Linguistic regression (T2F-LLR) variables to capture both randomness and higher-order fuzziness. 
We rigorously develop a solution framework based on a one-sigma confidence interval using the credibility measure to calculate the expected values and variances of the model output. 
To improve computational efficiency and usability of decisions, we introduce a heuristic algorithm tailored for non-meta datasets, significantly reducing the complexity of the model solving process. 
The experimental results demonstrate the effectiveness of our approach, which yields a mean absolute percentage error (MAPE) of weight equals 7.97\% with all variables statistically significant. We also provide the significant results for each product using the one-way analysis of variance test (one-way ANOVA test) ($p$-value = 0.15) and the paired $t$ test ($p$-value = 0.16). The results show that there is no significant difference between observed and predicted weights overall. 
This paper provides a robust and interpretable methodology for decision makers dealing with imprecise data and time-sensitive planning.
\end{abstract}


\IEEEpeerreviewmaketitle

\section{Introduction}\label{sect:1}

The task of making decisions is significantly influenced when uncertainty is present, encompassing both probabilistic and fuzzy measures. Probability measures trace their historical origins back to the 19th century, while the foundation of fuzzy measures lies within the realm of fuzzy sets. Real-life situations frequently involve the amalgamation of multiple uncertainties. To navigate this complex landscape, Baoding Liu introduced a credibility measure that seamlessly combines the domains of probability and fuzzy set approaches \cite{[A11]}. The concept of multi-uncertainty represented as a fuzzy random variable is credited to Kwakernaak \cite{[A12]} and \cite{[A13]}.

In this paper, we focus on harnessing expert knowledge to navigate the intricate terrain of multiple uncertainties. In the era of artificial intelligence, it becomes essential to adeptly process natural language expressions provided by individuals and effectively convey computational results in a manner that aligns with human communication.
The foundation of this study lies in the methodology of fuzzy control, which gained prominence in Japan during the 1980s. Leveraging this approach, we have endeavored to construct a natural language evaluation system that adeptly assesses structural damage caused by seismic events. In 1983, the powerful earthquake that struck San Francisco served as a significant example that had not received comprehensive scrutiny within the realm of multi-uncertainty. Such an undertaking falls within the purview of granular computing \cite{[1]}. Information granules are often formed based on resemblances, spatial or temporal adjacency, and functional attributes of constituent elements \cite{[4]}. 

The domain of fuzzy Linear regression (FLR) analysis has traditionally focused on datasets comprising numeric values, interval numbers, or fuzzy numbers, excluding stochastic elements. 
In real-world scenarios, the inherent indeterminacy associated with natural language has garnered substantial attention. Set theory is frequently employed to establish their values within specific cohorts or events. While it adeptly grapples with such uncertainties, it encounters challenges when confronted with these nuances.
The pioneering work of Lofti A. Zadeh in 1965 \cite{[8]} laid the foundation for the introduction of fuzzy sets, commonly known as Type-1 Fuzzy (T1F) sets. Subsequently, in 1982, H. Tanaka et al. introduced a linear programming (LP) approach to construct a fuzzy linear regression (FLR) model \cite{[7]}. Building upon these foundations, Watada and Tanaka extended the frontiers of the fuzzy quantification methodology.
A comprehensive overview of Watada's fuzzy models was presented in 2020. \cite{[A14]}.

It is noteworthy to acknowledge that the membership functions inherent in T1F sets have the capacity to embody uncertainty, as emphasized by Lofti A. Zadeh in 1975 through the introduction of Type-2 fuzzy (T2F) sets \cite{[10]}\cite{[10-1]}\cite{[10-2]}. However, despite their potential, the adoption of T2F sets has not gained widespread traction due to their intricate computational requirements and the challenges associated with comprehending and utilizing them. These limitations stem from two fundamental factors: firstly, the intrinsic three-dimensional nature of T2F sets, which renders their manipulation arduous, and secondly, the computational intricacies entailed in employing T2F sets compared to their Type-1 (T1) counterparts.

Before the 1980s, only a limited number of researchers, such as Mizumoto and Tanaka \cite{[11]}, delved into exploring the algebraic properties of different grades of T2F sets concerning negation, join, and meet operations. Dubois and Prade \cite{[12]} focused on investigating logical operations involving fuzzy-valued numbers. However, in recent times, T2F sets have gained increasing application in fuzzy logic systems to address both linguistic expressions and numerical fuzziness and randomness \cite{[13]}, \cite{[14]}, \cite{[15]}, \cite{[16]}, \cite{[17]}.

Numerous fuzzy models have been developed specifically to address qualitative data arising from fuzzy environments, often involving subjective estimations by human experts. Tanaka pioneered the fuzzy linear regression (FLR) model \cite{[7]}, while Tanaka and Watada, as well as Watada and Tanaka, introduced possibilistic linear regression through the possibility measure \cite{[20]}. 
Additionally, Chang ventured into the field of FLR with the least-squares criterion \cite{[22]}, employing fuzzy arithmetic to leverage weighted and least-squares criteria. Expanding further, Watada proposed models for T1F models, including linguistic terms in the FLR model. %\cite{[23]}.
The solution method was obtained as a heuristic approach \cite{[A14]}. A more theoretical approach to problems involving fuzzy numbers can be found in \cite{[A02]}. 

The exploration of T2F sets and their applications in linguistic and uncertain information processing has garnered significant attention in recent research. One foundational contribution is by Zhou et al. \cite{[T2F01]}, who introduced a type-1 (T1) OWA operator designed to aggregate uncertain information with fuzzy weights derived from linguistic quantifiers. This operator functions as an uncertain OWA, facilitating the aggregation of linguistic opinions and preferences in decision-making contexts involving linguistic weights, which is pertinent to the broader theme of handling linguistic uncertainty.
Building upon the fuzzy set extensions, Liu et al. \cite{[T2F02]} examined the cloud model, an innovative approach that integrates probability measures into fuzzy sets. This model captures both the randomness and fuzziness inherent in linguistic concepts by extending the degree of membership to a probabilistic framework, thus providing a more comprehensive representation of linguistic uncertainty. The comparative analysis between cloud models and other fuzzy set extensions underscores the significance of probabilistic considerations in linguistic modeling.

Further advancing the theoretical landscape, Niewiadomski \cite{[T2F03]} delved into the properties and formal aspects of T2F sets, particularly in the context of linguistic summarization of databases. The study emphasizes the role of discrete and continuous T2F sets in representing quantifiers, summarizers, and qualifiers, highlighting their capacity to encapsulate the nuanced uncertainty present in linguistic data summaries. This aligns with the overarching goal of effectively modeling linguistic information with inherent ambiguity.
In the realm of fuzzy clustering, Linda et al. \cite{[T2F04]} proposed a general T2F $C$-means (GT2 FCM) algorithm tailored for uncertain fuzzy clustering tasks. Extending the traditional IT2 FCM via the $\alpha$-planes representation theorem, this approach enhances the robustness of clustering in uncertain environments, which is essential for applications involving linguistic data characterized by ambiguity and imprecision.
Abdullah et al. \cite{[T2F05]} introduced a novel T2F set of linguistic variables specifically designed for the Fuzzy Analytic Hierarchy Process (FAHP). Their method was applied to safety and early warning assessments in hot and humid environments, demonstrating the practical utility of T2F sets in decision-making processes that involve linguistic uncertainty. The comparative results further validate the feasibility of employing T2F sets in such contexts.

Addressing regression problems under linguistic uncertainty, Song et al. \cite{[A00]} developed a type-2 (T2) linguistic random regression model grounded in credibility theory. This model constructs confidence intervals for fuzzy inputs and outputs, effectively capturing the hybrid uncertainty of primary and secondary fuzziness. The nonlinear programming framework facilitates linguistic assessments, illustrating the applicability of T2F sets in modeling complex linguistic relationships.
Sun et al. \cite{[T2F06]} extended the application of T2F sets to route evaluation for unmanned aerial vehicles, formulating the problem as a multi-criteria decision-making task under uncertainty. Their integrated approach leverages the diverse input types—numerical, probabilistic, and linguistic—highlighting the flexibility and effectiveness of T2F sets in handling complex, uncertain decision environments.
ALso, Boran et al. \cite{[T2F07]} provided a comprehensive review of linguistic summarization methods within the fuzzy set framework, emphasizing the evaluation techniques and current applications. The review underscores the importance of fuzzy set extensions, including T2F sets, in enhancing the expressiveness and interpretability of linguistic summaries, and points to open issues and future research directions.

In the context of data clustering, Comas et al. \cite{[T2F08]} proposed a method based on interval T2F predicates, enabling interpretable and efficient clustering of data with inherent uncertainty. Their approach demonstrates the capacity of interval T2F sets to extract meaningful knowledge from uncertain data, which is relevant for linguistic data analysis where ambiguity is prevalent.
Chen et al. \cite{[T2F09]} introduced a proportional interval T2 hesitant fuzzy TOPSIS method for linguistic decision-making under uncertainty. By encoding linguistic information as proportional interval T2 hesitant fuzzy sets, their approach effectively models individual comprehension and linguistic hesitation, providing a nuanced framework for decision-making processes that involve linguistic uncertainty.

In summary, these studies illustrate the evolving nature languages of T2F sets and their critical role in modeling, aggregating, and analyzing linguistic information under conditions of uncertainty. The integration of probabilistic elements, class attributes, and algorithms highlights the diversity and depth of research in this field. 
The concept of a T2F-LLR model that we propose integrates these concepts to enable the handling of complex, uncertain linguistic phenomena.

To accommodate linguistic terms, we propose a framework that includes vocabulary translation and vocabulary matching procedures. These procedures facilitate the transformation of linguistic terms into membership functions within $[0, 1]$, representing fuzzy expressions. Essentially, linguistic terms are translated into fuzzy numbers or sets, which are then integrated into the framework of approximate reasoning. Regression analysis on fuzzy models and/or fuzzy data, as cited in \cite{[7]}, 
In the research of \cite{[A02]} which is applied to navigate the mapping and assessment process conducted by experts. 
This process involves linguistic terms that capture the features and characteristics of the objectives, and resulting return to linguistic terms' expression encapsulating the overall assessment.

Motivated by the considerations mentioned above, this paper introduces a novel class of linguistic terms with T2F set of linguistic linear regression (T2F-LLR) model grounded in credibility theory. 
These models facilitate the handling of T2F random variables inputs and outputs. Our approach adopts the credibility theory originally proposed by Liu \cite{[24]} to precisely specify the variance and expected value for the T2F sets.
Subsequently, the expected value and variances of T2F random data is calculated and then forms into T2F-LLR model. 
Through comprehensive training, our model can effectively evaluate the input data and generate the corresponding output with T2F by using our heuristic algorithm.

The remainder of this work is structured as follows. 
Section \ref{Sec:3} provides the mathematical background necessary for constructing the T2F-LLR model. In Section \ref{Sec:4}, we give a conscientious process for building the T2F-LLR model along with the solution method to obtain the results for the model. After that, Section \ref{Sec:5} presents a numerical experiment for a real case study. Finally, Section \ref{Sec:6} concludes the paper, and some future studies are given as well.

\section{Mathematical Functions}\label{Sec:3}
In this section, we will introduce some basic mathematical functions that we will use for building our T2F-LLR model.
\subsection{Type-2 Fuzzy (T2F) Set}\label{sec:3.1}
%\subsubsection{T2F Set} 

Zadeh extended his fuzzy set theory by introducing Type-2 Fuzzy (T2F) Sets, as documented in \cite{[10]}. According to the findings from \cite{[26]}, T2F sets are defined as "sets where the degrees of membership are represented by Type-1 fuzzy sets (T1F sets)." 
We give a clear definition as follows.

\begin{definition}\label{Def:T2F}{\bf Type-2 Fuzzy (T2F) Set}

A T2F set, denoted as $\tilde{A}$, is a "second-order" fuzzy set with its membership function $\mu_{\tilde{A}}(x, \mu)$, where $x$ is the primary variable and $\mu$ is the secondary variable.  
For a universe of discourse $X$, a T2F set can be expressed as:

\begin{equation}\label{(2.1)}
\tilde{A}=\int_{x \in X}  \int_{u \in U_x} \frac{\mu_{\tilde{A}}(x,u)}{u}/x,
\end{equation}    
\end{definition}
where $U_{x} \subseteq [0,1]$ denote the primary membership grade set of the variable $x \in X$, where $u \in U_{x}$ $\forall x \in X$.

A Type-2 (T2) membership function $\mu_{\tilde{A}}(x,u)$, as elucidated in the references \cite{[27]}, \cite{[28]} and \cite{[28a]}, operates within the context of $\Re^3$. Nevertheless, it's crucial to acknowledge that representing and comprehending T2 membership functions pose substantial challenges in terms of visualization, manipulation, and interpretation, as discussed in reference \cite{[27]}.

Let us consider a universe of discourse $x$ to establish a mathematical definition for a T2F set. 
If $A\in X$ is a fuzzy set, then it is defined through $\mu_A: X \rightarrow [0,1]$, its membership function. Then, we can express a fuzzy set using the subsequent formulation:
\begin{equation}\label{(2.2)x}
A=\{(x,\mu_A (x)): \mu_A(x) \in [0.1], \forall x \in X \}
\end{equation}

It is important to emphasize that we have exact values for the membership grades of the set $A$. The collection of fuzzy sets belonging to the universe $X$ can be denoted as $\tilde{P}(X)$. Then, we notice that a fuzzy set where the membership grades also exhibit fuzziness, is a T2F set $\tilde{A}\in  X$.

The latter indicates that $\mu_{\tilde{A}}(x)$ represents a fuzzy set within the range of $[0,1]$ for all $x$, and this concept can be articulated as follows:
\begin{equation}\label{(2.3)}
\tilde{A} = \{(x,\mu_{\tilde{A}} (x)): \mu_{\tilde{A}} (x) \in \tilde{P}([0,1]), \forall x \in X \}
\end{equation}

This implies the following: 

For all $x$ in $X$, there is a set $U_{x} \subseteq [0,1]$ satisfying the membership function $\mu_{\tilde{A}}(x)$ that can be defined as $\mu_{\tilde{A}}(x): U_{x} \rightarrow [0,1]$. Applying equation (\ref{(2.2)x}), we obtain:
\begin{equation}\label{(2.4)}
\begin{array}{l}
\mu_{\tilde{A}} = \{(u,\mu_{\tilde{A}}(x) (u)) \\
\hspace{1cm} : \mu_{\tilde{A}} (x) (u) \in [0,1], \forall u \in U_{x} \subseteq [0,1]  \}
\end{array}
\end{equation}

We call $V_{x}$  the primary membership of $x$, and we refer to $\mu_{\tilde{A}}(x)$ as the secondary membership of $x$. Combining equations (\ref{(2.3)}) and (\ref{(2.4)}) yield the following expression:
\begin{equation}\label{(1.5)}
\begin{array}{l}
\tilde{A}=\{(x,(u,\mu_{\tilde{A}}(x)(u)))  \\
\hspace{1cm} : \mu_{\tilde{A}} (x) (u) \in [0,1], \forall x \in X \wedge  \forall u \in U_{x} \subseteq [0,1]  \}
\end{array}
\end{equation}

\subsection{Type-Reduction Representation}
The vertical representation of T2F sets is very important in defining the concept of embedded sets within the T2F sets. This representation serves as a fundamental tool for understanding the structure and properties of T2F sets, especially for determining the centroid, a key characteristic of a T2F set.
When dealing with any type of fuzzy sets, the conventional approach involves discretization as the initial step in creating a computer-based representation. Discretization requires the conversion of continuous sets into discrete sets through a slicing process. The rationale behind discretization arises from the inherent limitations of computers, which can only handle a finite number of slices, making direct processing of continuous fuzzy sets unfeasible. By discretizing fuzzy sets into manageable slices, computer-based systems become able to efficiently process and analyze them, enabling practical use.

With this concept, we can know that for T2F sets, a slice pertains to a plane falling into one of two distinct categories. The first category encompasses a plane intersecting the axis $x$ and running parallel to the plane $u$-$z$. The second category comprises a plane intersecting the axis $u$ and parallel to the plane $x$-$z$. These two plane types serve the purpose of extracting specific insights and characteristics from a T2F set, thereby facilitating a more profound comprehension of its behavior and attributes.

Hence, a vertical slice of a T2F set can be defined as follows:

\begin{definition}\label{Def:2}{\rm \bf Vertical Slice \cite{[27]}\cite{[37]}\cite{[21]}}

A vertical slice of a T2F set is a plane intersecting the axis $x$ and parallel to the plane $u-z$. 
\end{definition}

\begin{figure}[ht]
\centering
\includegraphics[width=3.5in, keepaspectratio]{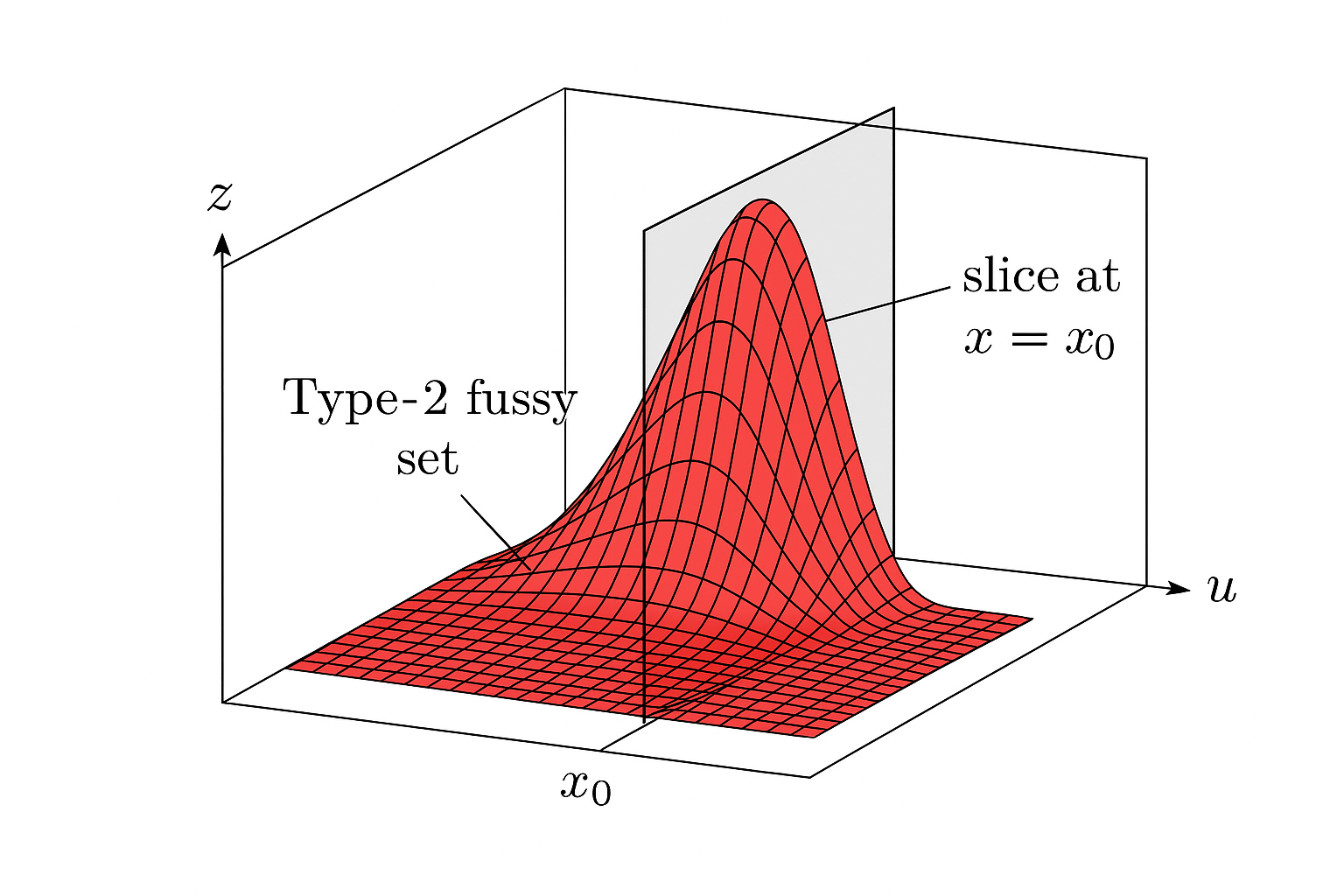}
\caption{Vertical Slice of T2F}\label{T2FVS}
\end{figure}

This particular type of slice enables an examination of membership grades at various $x$ values while maintaining a constant $u$-variable. The analysis of these vertical slices provides valuable insights into the fluctuating membership grades of the T2F set along axis $x$. We give this concept figure as shown in Figure \ref{T2FVS}.

After we have this concept, we can easily use mathematical definition to define the T2F set in the Euclide space and calculate each slice on the difference surface. Lucas \cite{{[37]}} employed this definition to build an easy way to deal with vague data and collected them by using a type reduction method which called vertical slice centroid type reduction (VSCTR). There were many research works, \cite{[37]}  \cite{[GRE00]} \cite{[36]} \cite{[39]}, explored the efficiency and accuracy of the VSCTR method in defuzififying the T2F data into T1F. Considering the exhaustive time, benchmark values, the method can efficiently approximates the doamin data. Hence, we will consider the method in dealing with our T2F data into our heuristic algorithm.

Let us introduce our type-reduce processes in the following subsection.

\subsection{Decomposing T2F sets into T1F Sets}

With previous concepts, we should have a rule to decompose the T2F sets into T1F sets. First, we introduce how to decompose the T2F set into discretization form.

\begin{definition}\label{Def:3}{\rm \bf Degree of Discretization}

A T2F set assigns to each element $x \in X$ a fuzzy set of membership grades (as show in equation (\ref{(1.5)})). This fuzzy set is usually represented over the interval $[0,1]$.
%— and that interval needs to be approximated in practice. 
The {\it Degree of Discretization} is the number of points or steps used to represent the secondary membership grades (i.e., the $z$-axis in 3D representation) for each value of $x$.

%The degree of discretization pertains to the separation between slices. 
In the case of a T2F set, the primary domain and the secondary domain undergo discretization. We divide the primary domain in vertical slices, while the secondary domain, spanning the $U = [0, 1]$, unit interval, can exhibit varying levels of discretization. 
It's essential to note that the degree of discretization in the secondary domain is not necessarily constant and can vary. 
Multiple strategies for discretizing T2F sets are discussed in \cite{[39]}.
\end{definition}

Furthermore, when addressing the process of transforming T1F sets into crisp values, various techniques can be applied. These techniques encompass methodologies such as the centroid calculation, the center of maxima determination, and the avearge of maxima computation, as discussed in \cite{[40]}. However, defuzzification of discretized T2F sets involves a dual-step procedure, as expressed in reference \cite{[41]}.

The initial step in this process is known as type-reduction, which serves to convert the T2F set into a T1F set. Subsequently, the second step involves the defuzzification of the resulting T1F set, employing the appropriate defuzzification technique tailored for T1F sets. This two-step approach empowers us to extract precise, non-fuzzy values from the discretized T2F sets, and analyze T2F sets in future decision-making scenarios. Furthermore, the type-reduction algorithm for T2F sets is based on the Extension Principle \cite{[10]}. The Extension Principle extends the operations originally defined for crisp numbers to T1F sets. Consequently, the defuzzification techniques for T2F sets are derived from and incorporate the methods used for type-1 (T1) defuzzification. This enables us to utilize established defuzzification methods for T1F sets when dealing with T2F sets, ensuring consistency and compatibility in the defuzzification process.

As we have shown the benefit of the VSCTR method in defuzififying the T2F data into T1F. We still need to know another concept for defuzzifying multi T2F sets into T1F sets.

An embedded T2F set (embedded set), also known as {\bf Wavy-Slice Representation Theorem} \cite{[27]}, has a very important role in the type-reduction process.  An embedded set establishes a unique link with the overarching T2F set as shown in the following: for each primary domain value $x$, a corresponding secondary domain value $u$ exists along with the secondary membership grade associated $\mu_{\tilde{A}}(x)(u)$, defined through both the primary domain and the secondary domain values. 

We give a comprehension mathematical expression for the Embedded Set in Definition \ref{Def:4}.
\begin{definition}\label{Def:4}{\rm\bf Embedded Set} \cite{[GRE01]}

Let $\tilde{A}$ be a T2F set belonging to the universe $X$. For the discrete universes of discourse $U$ and $X$, a type-2 (T2) embedded set $\tilde{A}_i$ of $\tilde{A}$ is defined as follows: 
\begin{equation}\label{(1.6)}
\begin{array}{ll}
\tilde{A}_i = & \{(x_i, (u_i, \mu_{\tilde{A}}(x_i)(u_i))) \\
& | \forall~ i \in \{1,\ldots,N\} : x_i \in X, u_i \in U_{x_i} \subseteq U \}, 
\end{array}
\end{equation}
where $\tilde{A}_i$ encompasses a solitary element drawn from $U_{x_1}, U_{x_2}, \ldots, U_{x_N}$, denoted as $u_1, u_2, \cdots, u_N$, each paired along their corresponding secondary grade, represented as $\mu_{\tilde{A}}(x_1)(u_1)$, $\mu_{\tilde{A}}(x_2)(u_2)$, $\ldots$, $\mu_{\tilde{A}}(x_N)(u_N)$. 
\end{definition}

In Mendel and John's research works \cite{[27]}, they demonstrated that T2F sets can be represented by the union of their own T2 embedded sets. The outcome result is usually referred to as the Wavy-Slice Representation Theorem or the T2F Set Representation Theorem. Notably, this theorem was derived independently of the Extension Principle \cite{[10]}, which is generalize crisp mathematical concepts to fuzzy sets, offering a conceptually simpler approach to handling T2F sets. By leveraging this theorem, we rewrite the Definition \ref{Def:4} based on the Extension Principle to be more straightforwardly, facilitating the manipulation, and analysis of T2F sets.

\begin{lemma}\label{lemma0} {\bf Embedded Set (Based on Extension Principle)}

Let $\tilde{A}_i^j$ represent its $j$th T2 embedded set for the T2F set $\tilde{A}$. From Definition \ref{Def:4}, the T2 embedded set can be written as follows: 
\begin{equation}
\begin{array}{ll}
\tilde{A}_i^j = & \{(x_i, (u_i^j, \mu_{\tilde{A}}(x_i)(u_i^j))) \\
& | \forall~ i \in \{1,\ldots,N\} : x_i \in X, u_i^j \in U_{x_i} \subseteq U \},
\end{array}
\end{equation}
where $\displaystyle{\tilde{A} = \sum_{j=1}^{n} \tilde{A}_i^j}$, $n$ is defined as $\displaystyle{n \equiv \prod_{i=1}^{N} M_i}$.
\end{lemma}

Note that, in the stage corresponding to the type reduction, the secondary membership grades for the T2F sets undergo a $t$-norm operation ($*$). It is important that in the case of T2F sets with general secondary membership functions, the product $t$-norm may not yield meaningful results. Therefore, it is advisable to avoid using the product $t$-norm in such cases. Instead, for the purposes of our research scenario, we employed the minimum $t$-norm as an alternative.
Therefore, it is crucial for the domain $X$ to have a numerical terms to ensure the meaningfulness of the type reduction stage. Consequently, the type reduction stage yields a T1F set within the universe of discourse $U$. In practical application, discretizing the secondary domain $U$ is also needed to compute the type reduction stage. The exhaustive type-reduction algorithm proves useful in computing the the type-reduction set for a T2F set. We will describe the detail in the following susection.

\subsection{Architecture of Type-Reduction Algorithm}\label{sec:3.2} 

The Representation Theorem of Mendel and John offers a clear and direct approach to T2 defuzzification \cite{[27]}. Although embedded sets may not be explicitly mentioned in Definition \ref{Def:4}, they can be implicitly present by Credibility Measure. When expressed as an algorithm, known as the exhaustive type-reduction algorithm, embedded sets are explicitly referenced. The term "exhaustive" signifies that this method processes each embedded set individually.
It means that each embedded set undergoes systematic defuzzification as a T1F set. Subsequently, the minimum secondary membership grade associated with that particular embedded set is merged together with the defuzzified value. By iterating this procedure for each embedded set, a series of ordered pairs is generated, composing the type-reduction set. It's important to highlight that while discretization introduces some level of approximation, the exhaustive method doesn't introduce any additional inaccuracies beyond the discretization phase.

Before we start introducing our type-reduction algorithm, we first give some important theorem and lemma for using in algorithm. Let's one-by-one step introduce in the following context.

\begin{theorem}{\bf Credibility Theory}\label{thm2}

Let us consider a universal Set $\Gamma$ and a possibility measure $Pos$, which operates over $P(\Gamma)$, the power set for $\Gamma$. Consider $Y: \Gamma \rightarrow \mathbb{R}$ is a function referred to a fuzzy variable over $\Gamma$ \cite{[47], [A10]}. Hence, $\mu_Y$, the possibility distribution  of $Y$, can be written as $\mu_Y(t) = Pos_Y(t)$, where $t \in \mathbb{R}$ represents the possibility of the event $\{Y = t\}$.

If we take $Y$ as a fuzzy variable with the possibility distribution $\mu_Y$, we can express the necessity, the credibility, and the possibility for the event $\{Y \geq r\}$ in the following manner:
\begin{equation}\label{eqn:11}
\begin{array}{ll}
Pos \{Y \geq r\} = & \displaystyle{\sup_{t \geq r} \mu_Y (t)}\\ 
Nec \{Y \geq r\} = & \displaystyle{1 - \sup_{t \leq r} \mu_Y (t)}\\
Cr \{Y \geq r\} = & \displaystyle{\frac{1}{2}} (1 - \sup_{t \leq r} \mu_Y (t) + \displaystyle{\sup_{t \geq r} \mu_Y (t)) }
\end{array}
\end{equation}
\end{theorem}

In Equation (\ref{eqn:11}), we calculate the credibility measure  as the mean value of the necessity and the possibility  measures, denoted as $Cr(\cdot) = \displaystyle{\frac{Nec(\cdot)+Pos(\cdot) }{2}}$. 
Note that, %as stated in \cite{23}, 
the credibility measure is a self-dual set function, which implies that, for any set $A\in P(\Gamma)$, we have the $Cr(A) = 1 - Cr(A^c)$, where $A^c$ means the supplement of $A$.

The integration of the credibility measure is motivated by the desire to establish a measure that combines the optimistic nature of possibility, reflecting a degree of overlap, and the pessimistic nature of necessity, indicating a level of inclusion. Using the credibility measure, the mean value of the fuzzy variable can be defined in the following manner:

\begin{definition}\label{Def:6}{\bf Expected Value of Fuzzy Variable}\\
The mean value of a fuzzy variable $X$ can be defined as follows:

\begin{equation}\label{eqn:12} 
E[X] = \int_0^{\infty} Cr\{X \geq r\}dr - \int_{-\infty}^0 Cr\{X \leq r\}dr,
\end{equation}
existed when each integral has a finite value.
\end{definition}

\begin{lemma}\label{lem1} 
Let us consider a triangular fuzzy variable denoted by $X = (c, a^l, a^r)_T$ with its membership function as follows: 

\begin{equation}\label{eqn:13} 
\mu_X (x) = \left\{ 
\begin{array}{ll} 
\displaystyle{\frac{x - a^l}{c- a^l}}, & ; a^l \leq  x \leq c \\
\displaystyle{\frac{a^r -x}{a^r-c}}, & ; c \leq  x \leq a^r\\
0, & ; otherwise.
\end{array}\right.
\end{equation}

Hence, using equation (\ref{eqn:11}) in Theorem \ref{thm2}, the possibility, the necessity, and the credibility measure of $X$ can be expressed as follows:

\begin{equation}\nonumber 
Pos\{X \geq r\}= \left\{ 
\begin{array}{ll} 
1 & ; r \leq c \\
\mu_X(r)=\displaystyle{\frac{a^r -r}{a^r-c}} & ; c \leq  r \leq a^r\\
0 & ; r>a^r.
\end{array}\right.
\end{equation}

\begin{equation}\nonumber 
Nec\{X \geq r\}= \left\{ 
\begin{array}{ll} 
1 & ; r \leq a^l \\
1-\mu_X(r)=\displaystyle{\frac{c-r}{c-a^l}} & ; a^l \leq  r \leq c\\
0 & ; r>c.
\end{array}\right.
\end{equation}
Hence, we have
\begin{equation}\nonumber 
Cr\{X \geq r\}= \left\{ 
\begin{array}{ll} 
1 & ; r \leq a^l \\
\displaystyle{\frac{1}{2}}[2-\mu_X(r)] & ; a^l \leq  r \leq c\\
\displaystyle{\frac{1}{2}}\mu_X(r) & ; c \leq r \leq a^r\\
0 & ; r>a^r
\end{array}\right.
\end{equation}

To reduce the formula, we get the formula of credibility measure for triangular fuzzy variable as follows:

\begin{equation}\nonumber 
Cr\{X \geq r\}= \left\{ 
\begin{array}{ll} 
1 & ; r \leq a^l \\
\displaystyle{\frac{1}{2}}(1+\displaystyle{\frac{c-r}{c-a^l}}) & ; a^l \leq  r \leq c\\
\displaystyle{\frac{1}{2}}(\displaystyle{\frac{a^r -r}{a^r-c}}) & ; c \leq r \leq a^r\\
0 & ; r>a^r
\end{array}\right.
\end{equation}

Hence, by Definition \ref{Def:6}, we obtain the expected value  
\begin{equation}\label{eqn:14} 
E[X] = \displaystyle{\frac{1}{4}(a^l+2c+a^r).}
\end{equation}
\end{lemma}

To continue, let us present the definition for a fuzzy random variable, as well as the definitions for the operators of the mean value and the variance. These definitions serve as the fundamental building blocks for deriving subsequent theoretical results related to fuzzy random variables.

We present the definition for a fuzzy random variable in the following manner:

\begin{definition}\label{Def:7}{\bf Fuzzy Random Variable} \cite{[A12], [A13], [47]}

Consider $(\Omega, \Sigma, Pr)$ to be a probability space, and let $F_v$ denote a collection of fuzzy variables over the possibility space $(\Gamma, P(\Gamma), Pos)$. We define a fuzzy random variable as the map $X: \Omega \rightarrow F_u$ satisfying the following property: the function $Pos_X(\omega) \in B$ is measurable with respect to $\omega$ for all Borel subsets $B$ in $\mathbb{R}$.
\end{definition}

Consider $X$ as a fuzzy random variable within $\Omega$, the probability space. Building upon the preceding definition, it becomes evident that $X(\omega)$ represents a fuzzy variable for every $\omega$ within $\Omega$. Additionally, $X$, a fuzzy random variable, can be classified as positive if the fuzzy variable $X(\omega)$ is positively oriented with an almost certain probability, for almost all $\omega$.

Let's give an example in the following:

\begin{example}\label{Example:2}

Let $V$ be a random variable within $(\Omega, \Sigma, Pr)$, the probability space. Let $X(\omega) = (U(\omega) + 2, U(\omega) - 2, U(\omega) + 6)_T$ for each $\omega \in \Omega$, where $X(\omega)$ represents a triangular fuzzy variable within a particular possibility space $(\Gamma, P(\Gamma), Pos)$.
\end{example}

Hence, the variable $X$ can be characterized as a triangular fuzzy random variable. Thus, it has been established that the mean value of $X(\omega)$, given by $E[X(\omega)]$, is a measurable function with respect to $\omega$  for any fuzzy random variable $X$ defined within $\Omega$, where $\omega$ belongs to the sample space \cite{[46]}. In other words, it can be treated as a random variable. Therefore, the mean value of the fuzzy random variable $X$, denoted as $E[X(\omega)]$, can be calculated using the credibility measure with respect to the probability space.

\begin{definition}\label{Def:8}{\bf  Expected Value by Credibility Measure} \cite{[46]}

Consider a fuzzy random variable $X$ within $(\Omega, \Sigma, Pr)$, a probability space. The mean value of $X$, denoted as $E[X]$, can be expressed in the following manner:
\begin{equation}\label{eqn:15}
\begin{array}{ll}
E[X] = & \displaystyle{\int_{\Omega} \Big[ \int_0^{\infty} Cr \{X(\omega) \geq r\} dr}\\
& \displaystyle{- \int_{-\infty}^0 Cr \{X(\omega) \leq r\} dr
\Big] Pr(d\omega)}
\end{array}
\end{equation}

Note that when the credibility measure is replaced with the probability measure, the above formula can be represented as follows:
\begin{equation}\label{eqn:16}
\begin{array}{ll}
E[X] & = \displaystyle{ \int_0^{\infty} Pr \{X(\omega) \geq r\} dr} \displaystyle{- \int_{-\infty}^0 Pr \{X(\omega) \leq r\} dr}\\
& = \displaystyle{ \int_{-\infty}^{\infty} X \mu_{X(\omega)}d{\omega}}
\end{array}
\end{equation}
\end{definition}
%{\noindent under the condition of the existence of the integrals.}

\begin{lemma}\label{lem2}
Suppose that we have triangular fuzzy variables $X_i=(c_i, a^l_i, a^r_i)_T$ with probability $p_i$, such that $ \sum _{i=1}^{n} p_i=1, \forall~ i=1,2,\dots,n$. From Lemma \ref{lem1}, we know that $E(X_i)= \displaystyle{\frac{1}{4}(a^l_i+2c_i+a^r_i)}$. Hence, we can obtain the expected value by credibility measure from Definition \ref{Def:8} as follows:
\begin{equation}\label{eqn:17}
\begin{array}{ll}
E[X] = \displaystyle{\sum _{i=1}^n p_i E(X_i)}
= \displaystyle{\sum _{i=1}^n p_i \displaystyle{\frac{1}{4}(a^l_i+2c_i+a^r_i).}}
\end{array}
\end{equation}
\end{lemma}

Let us show an example to calculate $E(X)$ with triangular fuzzy variables.

\begin{example}\label{Example:3}
Let us consider a discrete random variable $V \in \Omega$, with values of $V_1 = 3$ having a probability of $0.2$, and $V_2 = 6$ having a probability of 0.8. Now, let us determine the mean value $E[X]$ by credibility measure. 

From Example \ref{Example:2}, we can derive $X(V_1) = (5, 1, 9)_T$ having a probability of $0.2$, and $X(V_2) = (8, 4, 12)_T$ having a probability of $0.8$. 
Hence, by equation (\ref{eqn:14}), we have:

\begin{equation}\nonumber 
\begin{array}{ll} 
E[X(V_1)] = & \displaystyle{\frac{(1 + 9 + 2 \times 5)}{4}} = 5, \\ \\
E[X(V_2)] = & \displaystyle{\frac{(4+ 12 + 2 \times 8)}{4}} = 8.
\end{array}
\end{equation}

Hence, from equation (\ref{eqn:17}), we have: 
\begin{equation}\nonumber 
\begin{array}{ll} 
E[X]= & 0.2 \times   E[X(V_1)] + 0.8 \times E[X(V_2)] = 7.4. 
\end{array}
\end{equation}
\end{example}

\begin{definition}\label{Def:8-1}{\bf Variance by Credibility Measure} \cite{[46]}

Consider a fuzzy random variable $X$ within $(\Omega, \Sigma, Pr)$, a probability space. The expected value by credibility measure is $E(X)$ as shown in Definition \ref{Def:8}. The variance of $X$, denoted as $Var[X]$, can be defined as follows: 
\begin{equation}\label{eqn:18}
Var(X) = E[(X-E(X))^2]
\end{equation}
Suppose that $Y=(X-E(X))^2$, to find the formula of $Var(X)$ is equal to find the $E(Y^2)$.
Hence, from equation (\ref{eqn:15}), we obtain the following formula:
\begin{equation}\label{eqn:19}
\begin{array}{ll}
E[Y^2] = & \displaystyle{\int_{\Omega} \Big[ \int_0^{\infty} Cr \{Y^2(\omega) \geq r\} dr}\\
& \displaystyle{- \int_{-\infty}^0 Cr \{Y^2(\omega) \leq r\} dr
\Big] Pr(d\omega)}
\end{array}
\end{equation}
Replace the credibility measure to the probability measure as shown in equation (\ref{eqn:16}), we obtain the following formula:
\begin{equation}\label{eqn:20}
\begin{array}{ll}
E[Y^2] & = \displaystyle{ \int_0^{\infty} Pr \{Y^2(\omega) \geq r\} dr} \displaystyle{- \int_{-\infty}^0 Pr \{Y^2(\omega) \leq r\} dr}\\
& = \displaystyle{ \int_{-\infty}^{\infty} Y^2 \mu_{Y^2(\omega)}d{\omega}}
\end{array}
\end{equation}
\end{definition}

\begin{lemma}\label{lem3}

Suppose that $E(X)=m= \displaystyle{\sum _{i=1}^n p_i \displaystyle{\frac{1}{4}(a^l_i+2c_i+a^r_i)}}$ from equation (\ref{eqn:17}), and

\begin{equation}\label{eqn:21}
\begin{array}{ll}
Y_i & =X_i-E(X)\\
& =(c_i-m, a^l_i-m, a^r_i-m)_T 
= (\hat{c}_i, \hat{a}^l_i, \hat{a}^r_i) 
\end{array}
\end{equation}
is a new triangular fuzzy variables with probability $p_i$, such that $ \sum _{i=1}^{n} p_i=1, \forall~ i=1,2,\dots,n$. 
Based on the same concept to produce the variance by credibility measure for the triangular fuzzy random variable, we obtain the following formula:
\begin{equation}\label{eqn:22}
\begin{array}{ll}
Var[X] & = \displaystyle{\sum _{i=1}^n p_i Var(X_i)}
= \displaystyle{\sum _{i=1}^n p_i E[(X_i-E(X))^2]}\\
 & = \displaystyle{\sum _{i=1}^n p_i E(Y_i^2)}
 = \displaystyle{\sum _{i=1}^n p_i E(Z_i),}
\end{array}
\end{equation}
\end{lemma}
where $Z_i=Y_i^2$, $Y_i$ is still a triangular fuzzy variable,
and $Z_i$ is no longer a triangular number (unless $Y_i$ is strictly non-negative or non-positive), but it still is a fuzzy number. 
Hence, we can calculate the possibility distribution $\mu_{Y_i^2}$ from Theorem \ref{thm2}. Furthermore, we can obtain the $Cr\{Y_i^2 \geq r\}$ and $E(Y^2)$ by using the same processes as shown in Definition {\ref{Def:6}} and Lemma {\ref{lem1}}. 

Let us show how to do it one by one step.

{\bf Step 1: Find Possibility distribution $\mu_{Y^2}(z)$}

We have $Z=Y^2$, where $Y= (\hat{c}, \hat{a}^l, \hat{a}^r)$ as defined in equation (\ref{eqn:21}). To simply present the formula, we deleted the $i$ symbol here. By Lemma \ref{lem1}, we have membership function of $Y$ as follows:
\begin{equation}\label{eqn:23}
\mu_Y (y) = \left\{ 
\begin{array}{ll} 
\displaystyle{\frac{y - \hat{a}^l}{\hat{c}- \hat{a}^l}} & ; \hat{a}^l\leq  y \leq \hat{c} \\
\displaystyle{\frac{\hat{a}^r-y}{\hat{a}^r-\hat{c}}} & ; 
\hat{c} \leq  y \leq \hat{a}^r\\
0 & ; otherwise.
\end{array}\right.
\end{equation}
Sine $z=y^2$, it reduce as $y=\pm \sqrt{z}$.
We have 
\begin{equation} \label{eqn:24}
\begin{array}{ll}
\mu_{Y^2} (z) &= sup \{\mu_Y (y) | y \in \Re, y^2 =z\} \\ 
& = \max \{ \mu_Y (\sqrt{z}), \mu_Y (-\sqrt{z}) \}
\end{array}
\end{equation}
This is a general inverse-mapping method. But since $\mu_Y(y)$ is piecewise function, we can explicitly calculate it for each root.

{\bf Case 1:} $\sqrt{z} \in [\hat{a}^l, \hat{c}]$\\
Then, 
\begin{equation} \nonumber 
\mu_Y (\sqrt{z})= \displaystyle{\frac{\sqrt{z} - \hat{a}^l}{\hat{c}- \hat{a}^l}}
\end{equation}

{\bf Case 2:} $\sqrt{z} \in [\hat{c}, \hat{a}^r]$\\
Then, 
\begin{equation} \nonumber 
\mu_Y (\sqrt{z})= \displaystyle{\frac{\hat{a}^r-\sqrt{z}}{\hat{a}^r-\hat{c}}}
\end{equation}

{\bf Case 3:} $-\sqrt{z} \in [\hat{a}^l, \hat{c}] \Rightarrow \sqrt{z} \in [-\hat{c}, -\hat{a}^l]$\\
Then,
\begin{equation} \nonumber 
\mu_Y (-\sqrt{z})= \displaystyle{\frac{-\sqrt{z} - \hat{a}^l}{\hat{c}- \hat{a}^l}}
\end{equation}

{\bf Case 4:} $-\sqrt{z} \in [\hat{c}, \hat{a}^r] \Rightarrow \sqrt{z} \in [-\hat{a}^r, -\hat{c}]$\\
Then, 
\begin{equation} \nonumber 
\mu_Y (-\sqrt{z})= \displaystyle{\frac{\hat{a}^r+\sqrt{z}}{\hat{a}^r-\hat{c}}}
\end{equation}

Hence, when $z>0$, we have
\begin{equation} \label{eqn:25}
\begin{array}{ll}
\mu_{Y^2} (z) & = \max \{ \mu_Y (\sqrt{z}), \mu_Y (-\sqrt{z}) \} \\
 & = \max \left\{ \left\{
\begin{array}{ll}
\displaystyle{\frac{\sqrt{z} - \hat{a}^l}{\hat{c}- \hat{a}^l}} & ; \hat{a}^l\leq \sqrt{z} \leq \hat{c} \\
\displaystyle{\frac{\hat{a}^r-\sqrt{z}}{\hat{a}^r-\hat{c}}} & ; 
\hat{c} \leq  \sqrt{z} \leq \hat{a}^r\\
 0 & ; otherwise.
\end{array}\right.\right.,~\\
& ~~~\left\{
\begin{array}{ll}
\displaystyle{\frac{-\sqrt{z} - \hat{a}^l}{\hat{c}- \hat{a}^l}} & ; \hat{a}^l\leq -\sqrt{z} \leq \hat{c} \\
\displaystyle{\frac{\hat{a}^r+\sqrt{z}}{\hat{a}^r-\hat{c}}} & ; 
\hat{c} \leq  -\sqrt{z} \leq \hat{a}^r\\
 0 & ; otherwise.
\end{array}\right\}.
\end{array}
\end{equation}

{\bf Step 2: Find Credibility Measure $Cr\{Y^2 \geq r\}$}

Suppose that we get the membership function $\mu_{Y^2}(z)$ as follows: (For the case $0<\hat{a}^l<\hat{c}<\hat{a}^r$)
\begin{equation}\nonumber 
\mu_{Y^2} (z)  = \left\{
\begin{array}{ll}
\displaystyle{\frac{\sqrt{z} - \hat{a}^l}{\hat{c}- \hat{a}^l}} & ; \hat{a}^l\leq \sqrt{z} \leq \hat{c} \\
\displaystyle{\frac{\hat{a}^r-\sqrt{z}}{\hat{a}^r-\hat{c}}} & ; 
\hat{c} \leq  \sqrt{z} \leq \hat{a}^r\\
 0 & ; otherwise.
\end{array}\right.
\end{equation}

Note that, although the distribution $\mu_{Y^2}(z)$ is not an triangular fuzzy variables, it still can be presented as piecewise function. We can see that, in the first period, the function  $\displaystyle{\frac{\sqrt{z} - \hat{a}^l}{\hat{c}- \hat{a}^l}}$ is a monotonic increment function during the periods $(\hat{a}^l)^2 \leq z \leq \hat{c}^2$, and the function $\displaystyle{\frac{\hat{a}^r-\sqrt{z}}{\hat{a}^r-\hat{c}}}$ is a monotonic descending function during the periods $\hat{c}^2 \leq z \leq (\hat{a}^r)^2$. It keeps the same rule as we shown in Lemma \ref{lem1} that implied the credibility measure as follows:
\begin{equation}\label{eqn:26}
Cr\{Y^2 \geq r\}= \left\{ 
\begin{array}{ll} 
1 & ; r \leq (\hat{a}^l)^2 \\
\displaystyle{\frac{1}{2}}[2-\mu_{Y^2}(r)] & ; (\hat{a}^l)^2 \leq  r \leq \hat{c}^2\\
\displaystyle{\frac{1}{2}}\mu_{Y^2}(r) & ; \hat{c}^2 \leq r \leq (\hat{a}^r)^2\\
0 & ; r>(\hat{a}^r)^2
\end{array}\right.
\end{equation}

Note that, in the domain which is in negative value (for $\hat{a}^l<\hat{c}<\hat{a}^r<0$), the square roots period will flip orientation. There are some special cases such as $0 \in [\hat{a}^l, \hat{c}]$ or $0 \in [\hat{c}, \hat{a}^r]$, you must carefully check the domain of both $\pm \sqrt{z}$ in each cases.

{\bf Step 3: Find $E(Y^2)$ and $Var(X)$}
Base on the Definition \ref{Def:8} and equation \ref{eqn:19}. Hence, we can get the $Var(X)$.

We give an example for a special case (including 0 in the domain).

\begin{example}\label{Example:4}%  4 
 As the same data in Example {\ref{Example:3}}, we have 
calculated the $E(X)=7.4$. Now, we want to calculate the $Var(X)$ by credibility measure. It means that we need to calculate $Var(X)=E[(Y_1)^2] \times 0.2+ E[(Y_2)^2] \times 0.8$, where $Y_i=X_i-E(X)=X_i - 7.4, \forall ~i=1,2$. 

From Example {\ref{Example:3}}, we have:
\begin{equation}\nonumber 
\begin{array}{ll} 
Y_1 = & (5-7.4, 1-7.4, 9-7.4)_T = (-2.4, -6.4, 1.6)_T\\
Y_2 = & (8-7.4, 4-7.4, 12-7.4)_T = (0.6, -3.4, 4.6)_T
\end{array}
\end{equation}

To calculate the $E[(Y_i)^2], ~ \forall~ i=1,2$, we first need to calculate the membership function $\mu_{Y_i^2}(x), ~ \forall~ i=1,2$.

As we show in equation (\ref{eqn:24}), we need to find out
\begin{equation}\nonumber 
\begin{array}{ll} 
\mu_{Y^2_1}(x) = & \max \{ \mu_{Y_1}(\sqrt{x}), \mu_{Y_1}(-\sqrt{x})\}\\
\mu_{Y^2_2}(x) = & \max \{ \mu_{Y_2}(\sqrt{x}), \mu_{Y_2}(-\sqrt{x})\}
\end{array}
\end{equation}

Let us first calculate $\mu_{Y_1^2}(x)$.

Since $Y_1 = (-2.4, -6.4, 1.6)_T$, we have $-\sqrt{x} \in [-6.4, -2.4]$ (Case 3). Hence, we have left-hand membership function as follows:
\begin{equation} \nonumber
\begin{array}{ll} 
\mu_{Y_1} (-\sqrt{x})= & \displaystyle{\frac{-\sqrt{x} - (-6.4)}{-2.4 -(-6.4)}}= \displaystyle{\frac{6.4-\sqrt{x}}{4}}\\
& for~ -6.4 \leq -\sqrt{x} \leq -2.4.
\end{array}
\end{equation}
it implied to
\begin{equation} \nonumber
\begin{array}{ll} 
\mu_{Y_1} (-\sqrt{x})= & \displaystyle{\frac{6.4-\sqrt{x}}{4}}, ~~for~~(2.4)^2 \leq x \leq (6.4)^2.
\end{array}
\end{equation}

Second, we check the right-hand membership function. Since $0 \in [-2.4, 1.6]$, we need to find two intervals: $[-2.4, 0]$ and $[0, 1.6]$. 
$-\sqrt{x} \in [-2.4, 0]\in [-2.4, 1.6]$ (Case 4). Hence, we have a right-hand membership function as follows:
\begin{equation} \nonumber
\begin{array}{ll}
\mu_{Y_1} (-\sqrt{x})= & \displaystyle{\frac{1.6+\sqrt{x}}{1.6-(-2.4)}} = \displaystyle{\frac{1.6+\sqrt{x}}{4}}\\
& for~ -2.4 \leq -\sqrt{x} \leq 0.
\end{array}
\end{equation}
it implied to
\begin{equation} \nonumber
\begin{array}{ll} 
\mu_{Y_1} (-\sqrt{x})=  \displaystyle{\frac{1.6+\sqrt{x}}{4}}, ~~ for ~~ 0 \leq x \leq (2.4)^2.
\end{array}
\end{equation}

We still have another interval that is $\sqrt{x} \in [0, 1.6]\in [-2.4, 1.6]$ (Case 2). Hence, we have a right-hand membership function as follows:
\begin{equation} \nonumber
\begin{array}{ll}
\mu_{Y_1} (\sqrt{x})= & \displaystyle{\frac{1.6-\sqrt{x}}{1.6-(-2.4)}} = \displaystyle{\frac{1.6-\sqrt{x}}{4}}\\
& for~ 0 \leq \sqrt{x} \leq 1.6.
\end{array}
\end{equation}
it implied to
\begin{equation} \nonumber
\begin{array}{ll} 
\mu_{Y_1} (\sqrt{x})=  \displaystyle{\frac{1.6-\sqrt{x}}{4}}, ~~ for ~~ 0 \leq x \leq (1.6)^2.
\end{array}
\end{equation}

Since we should choose the maximum $\max \{ \mu_{Y_1}(\sqrt{x}), \mu_{Y_1}(-\sqrt{x})\}$ for $\mu_{Y^2_1}(x)$.
Hence, we have 
\begin{equation} \nonumber
\begin{array}{ll} 
\mu_{Y^2_1}(x)= & \max \{ \displaystyle{\displaystyle{\frac{1.6+\sqrt{x}}{4}}, \frac{1.6-\sqrt{x}}{4}}\}\\
& = \displaystyle{\frac{1.6+\sqrt{x}}{4}},~~for ~~ 0 \leq x \leq (2.4)^2.
\end{array}
\end{equation}

Finally, we get the membership function $\mu_{Y^2_1}(x)$ as follows: 
\begin{equation} \nonumber 
\mu_{Y_1^2} (x) = \left\{ \begin{array}{ll}
     \displaystyle{\frac{(1.6 + \sqrt{x}}{4}} &;  0 \leq t < 2.4^2 \\  \\
     \displaystyle{\frac{(6.4 - \sqrt{x})}{4}} &;  2.4^2 \leq t < 6.4^2 \\ \\ 
      0                               &; otherwise.  
\end{array} \right.
\end{equation}
\end{example}

Furthermore, we calculate 
\begin{equation} \nonumber
Cr\{Y_1^2 \geq r\}  = \left\{ \begin{array}{ll}
     \displaystyle{\frac{(2- \mu_{Y_1^2} (r))}{2}} &;  0 \leq t \leq 2.4^2 \\ \\ 
     \displaystyle{\frac{\mu_{Y_1^2} (r)}{2}} &;  2.4^2 \leq t \leq 6.4^2 \\ \\ 
      0                               &; otherwise.  
\end{array} \right.
\end{equation}

Therefore, from Definition \ref{Def:6}, we found $E[(X(V_1)-7.4)^2]= E[Y_1^2]$ in the following manner: 

\begin{equation} \nonumber 
 \begin{array}{ll}
E[Y_1^2] &= \displaystyle{\int_0^{\infty} Cr \{Y_1^2 \geq r\}} dr \\ \\
         &= \displaystyle{\int_0^{2.4^2} \frac{1}{2}(2- \frac{1.6+\sqrt{t}}{4})} dr
           + \displaystyle{\int_{2.4^2}^{6.4^2} \frac{1}{2}(\frac{6.4-\sqrt{t}}{4})}dr \\ \\
        &= 12.08  
\end{array}  
\end{equation}

Likewise, we can get $E[(X(V_2) - 7.4)^2]= E[Y_2^2]= 4.25$. Hence, $Var(X)= 0.2 \times   E(Y_1^2) + 0.8 \times  E(Y_2^2)^2 
= 0.2 \times 12.08 + 0.8 \times 4.25 = 5.81.$

\section{Building T2F-LLR Model under Multi Uncertainty}\label{Sec:4}

A T2F-LLR model is built by considering linguistic terms in a linear regression (LR) model. To describe the linguistic terms in a numerical equation, we give linguistic terms in each different case as having a membership function with respect to some weight. Hence, T2F random variables are considered in the problem. 
In our scenario problem, we would like to find a goal function that could provide the best investment strategy and give the decision-maker a guideline in advance. To solve this problem, we set up a T2F-LLR model under Multi Uncertainty. 

In the following subsections, we explain (1) how to transfer the linguistic terms into T2F random data and illustrate (2) how to solve the T2F-LLR model , then (3) we give the criteria rule for the heuristic algorithm. 

\subsection{T2F-LLR model (Credibility Based with Confidence Interval)}\label{sec:3.6}
Suppose that all linguistic terms have been transformed into fuzzy random data. The goal function in our scenario is generated into a T2F-LLR model, which is specifically designed for triangular fuzzy random data with probability in our case.
The parameters in the T2F-LLR model are shown in Table \ref{table2-1}, where $X_{ij}$ represents the input data with triangular fuzzy random data that are derived from linguistic terms, and $Y_i$ represents the output data which also are in triangular fuzzy data, for all $i  \in\{1,2,\dots,N\}$ and $j  \in\{1,2,\dots, J\} $. We give clearer definition as follows:

\begin{definition}\label{Def:9}{\rm\bf {Fuzzy Random Variables with Probability}}\\
Suppose that $X_{ij}$ and  $Y_i$ are all triangular fuzzy random data from linguistic terms, the parameters of Fuzzy Random Variables with Probability are defined as follows:
\begin{equation}
\begin{array}{ll}
X_{ij}  = & \bigcup_{t=1}^{M_{X_{ij}}} \{(X_{ij}^t,X_{ij}^{tl},X_{ij}^{tr}), q_{ij}^t \},  \\
Y_i = & \bigcup_{t=1}^{M_{Y_{i}}} \{(Y_i^t,Y_i^{tl},Y_i^{tr}), p_i^{t} \},
\end{array}
\end{equation}
where $(X_{ij}^t,X_{ij}^{tl},X_{ij}^{tr})$ and $(Y_i^t,Y_i^{tl},Y_i^{tr})$ are corresponding with the probabilities $q_{ij}^t$ and $p_i^{t}$, for all $i = 1,2,\dots,N $, $j= 1,2,\dots,J $, $t = 1,2,\dots,M_{X_{ij}}$, and $t = 1,2,\dots,M_{Y_i}$, respectively.
\end{definition}

\begin{table}[ht]
\caption{Linquistic terms for Each Sample $i$ provided by Experts}\label{table2-1}
\centering
\begin{tabular}{c|ccc|c}
\hline
sample & \multicolumn{3}{c|}{Input Value with Probability} & Output Value \\  
$i$ & \multicolumn{3}{c|}{$(X_{ij}, p_{ij})$}&  $Y_i$ \\  
\hline
1 & $(X_{11}, p_{11})$ & $\cdots$ & $(X_{1J}, p_{1J})$ &  $(Y_{1}, p_{1})$  \\
2 & $(X_{21}, p_{21})$ & $\cdots$ & $(X_{2J}, p_{2J})$ &  $(Y_{2}, p_{2})$  \\
3 & $(X_{31}, p_{31})$ & $\cdots$ & $(X_{3J}, p_{3J})$ &  $(Y_{3}, p_{3})$  \\
%4 & (good,0.2) & $\cdots$ & (very good,0.1) & (good,0.1) \\ 
$\vdots$ & $\vdots$ & & $\vdots$ &$\vdots$  \\
$i$ & $(X_{i 1}, p_{i 1})$ & $\cdots$ & $(X_{i j}, p_{i j})$ &  $(X_{i}, p_{i})$ \\ 
$\vdots$ & $\vdots$ & &$\vdots$ & $\vdots$  \\
N & $(X_{N1}, p_{N1})$ & $\cdots$ & $(X_{NJ}, p_{NJ})$ &  $(Y_{N}, p_{N})$  \\ \hline
\end{tabular}
\end{table}

Now, we have the T2F-LLR model by using fuzzy coefficients $\bar{A}_1,\cdots, \bar{A}_J$ as follow formula:

\begin{equation}\label{(2.5)}
\bar{Y}_i = \bar{A}_1 X_{i1} +\bar{A}_2 X_{i2} +\cdots+\bar{A}_J X_{iJ} 
\end{equation}

To define a confidence interval for a fuzzy random variable, we rely on the mean value ($\mu_X$) and variance ($\sigma^2_X$) of the fuzzy random variable, utilizing principles from credibility theory. If we consider a one-sigma confidence interval, it can be represented in the following definition.

\begin{definition}\label{Def:10}{\rm\bf {One-Sigma Confidence Interval}}\\
Consider a fuzzy random variable $X$ with mean value $\mu_X$ and variance $\sigma^2_X$. The one-sigma confidence interval for $X$ is defined as:
\begin{center}
$[\mu_X - \sigma_X, \mu_X + \sigma_X]$
\end{center}
This interval represents a range centered around the expected value $\mu_X$, within which $X$ is expected to lie with a certain level of confidence.
\end{definition}

In equation (\ref{(2.5)}), the confidence intervals are utilized to capture the uncertainty associated with the fuzzy parameters $\bar{A}_i$ in the T2F-LLR model, which considers both the fitness (how well it fits the data) and the fuzziness caught by the T2F-LLR model.

Various approaches can be employed to establish the inclusion $\supset_h$ operation for fuzzy random variables described as confidence interval. One common approach is to combine the mean value and variance of fuzzy random variables and ensure that the inclusion operation is valid for a specified confidence level $h$. This facilitates the representation of uncertainty and the capture of the range of possible values for fuzzy parameters.
It is important to note that there are other methods for establishing the fuzzy random inclusion operation $\supset_h$, which can result in other intricate FRLR models. In our T2F-LLR model, we strive to preserve more comprehensive information about fuzzy random data and directly consider the fuzzy inclusion operation in the product between a fuzzy value and a fuzzy parameter at a specific level of probability.

By incorporating the confidence interval and defining the fuzzy random inclusion operation for the T2F-LLR model, we provide a comprehensive framework for modeling and analyzing data under uncertain and fuzzy conditions.

In our scenario, we consider triangular fuzzy random number in T2F-LLR model, we consider the One-Sigma Confidence Interval based on Definition \ref{Def:10} and the expected value and variance in Definition \ref{Def:8} and \ref{Def:8-1}
In the case of triangular fuzzy random data for the type 2 fuzzy numbers the following formula can be calculated by referring to Lemmas \ref{lem2} and \ref{lem3}. 

\begin{equation}\label{(2.6)}
I(e_{X_{ij}}, \sigma_{X_{ij}}) = [E(X_{ij})-\sqrt{Var(X_{iJ})},E(X_{iJ})+\sqrt{Var(X_{iJ})}] 
\end{equation}

For obtaining natural language results, we apply the matching process to the obtained fuzzy numbers in order to select the fittest natural language expression (Vocabulary Matching). To begin, let us examine the one-sigma confidence interval for any fuzzy random variable in the following manner: 
\begin{equation}\label{(2.7)}
\begin{array}{rl}
I(e_{X_{iJ}} , \sigma_{X_{iJ}})= & [e_{X_{iJ}}- \sigma_{X_{iJ}}, e_{X_{iJ}} + \sigma_{X_{iJ}} ]\\ 
I(e_{Y_i} , \sigma_{Y_i})= & [e_{Y_i} - \sigma_{Y_i} , e_{Y_i} + \sigma_{Y_i}] 
\end{array}
\end{equation}

Next, we construct the fuzzy random linear regression model based on confidence intervals in the following manner: 
\begin{equation} \label{(2.8)}
\left. \begin{array}{rl}
\min_{\bar{\textbf{A}}} & J(\bar{\textbf{A}}) =  \sum_{j=1}^j (\bar{\textbf{A}}_j^r -\bar{\textbf{A}}_j^l ) \\ 
~& ~\\
\text{subject to } & \bar{\textbf{A}}_j^r \geq \bar{\textbf{A}}_j^l \\ 
& \bar{Y}_i = \sum_{j=1}^J \bar{\textbf{A}}_j I(e_{X_{iJ}}, \sigma_{X_{iJ}}) \supset_h I(e_{Y_i}, \sigma_{Y_i}) \\ 
%(2.8)
\end{array}\right\}
\end{equation}
where $i  \in\{1,2,\dots,N\} $, $j  \in\{1,2,\dots,J\} $, while the fuzzy inclusion relation established at level $h$ is denoted by $\supset_h$ . As the product of a confidence interval and a fuzzy number (or coefficient) is decided by the signs of every component, solving the model (\ref{(2.8)}) involves considering all cases that correspond to the to various combinations of signs for the $\sigma$-confidence intervals of the fuzzy random data and the fuzzy coefficients.

\subsection{Solution of the T2F-LLR model}\label{sec:3.7}
Upon expressing input $X_{ij}$ through its left and right endpoints from the expected primary grade intervals $\underline{e}_{ij}$ and $\bar{e}_{ij}$, respectively, the T2F linear regression model can be obtained in the model to incorporate the average interval values across every sample. Hence, it is both necessary and sufficient to only take into consideration the two vertices corresponding to the endpoints of the interval for all dimension for a given sample to solve the model. As an example, a single sample characterized by a single input interval feature can be represented by two vertices, which correspond to the endpoints of the interval, in addition to a fuzzy output value.

In the proposed model (\ref{(2.8)}), it is possible to reformulated the solution through a problem involving N samples, a single output, and K input interval values. However, solving such problem poses challenges, primarily arising from the $NK$ products involving fuzzy coefficients and confidence intervals. To address this problem, one can utilize a vertices method. By introducing multidimensional vertices we can create new sample points with fuzzy output values. Through this technique, conventional methods can be used to effectively find a solution for this problem.

However, it is important to note that this problem faces a combinatorial explosion, which becomes more prominent as the quantity of variables rises. The complexity of this problem grows exponentially, making it computationally demanding.

By considering the mean interval values and utilizing the vertices method, the T2F linear regression model allows for an effective representation of the samples and reduces the computational burden by focusing on the endpoints of the intervals in each dimension.
\begin{equation}\label{(2.9)}
\begin{array}{ll}
 \bar{e}_{ij} & = E(X_{ij}) + \sqrt{Var(X_{ij})}=e_{X_{ij}}+ \sigma_{X_{ij}},\\
 \underline{e}_{ij} & = E(X_{ij}) - \sqrt{Var(X_{ij})}=e_{X_{ij}}- \sigma_{X_{ij}}
\end{array}
\end{equation}
or $i  \in\{1,2,\dots,N\} $; $j  \in\{1,2,\dots,J\} $; the initial model can be transformed to a conventional fuzzy linear regression model through the application of a vertices method as follows:

\begin{equation}\label{eq:2.10}
\scriptstyle
\left.
\begin{array}{rl}
\displaystyle{\min_{\bar{A}}}  
& \displaystyle{J(\bar{A}) =  \sum_{j=1}^J (\bar{A}_j^r - \bar{A}_j^l )} \\[6pt] 
\text{subject to} & \bar{A}_j^r \geq \bar{A}_j^l \\[6pt] 
(1) \rightarrow & \bar{Y}_i = \bar{A}_1 \underline{e}_{i1} + \bar{A}_2 \underline{e}_{i2} + \cdots + \bar{A}_J \underline{e}_{iJ} \\
& \qquad \qquad \qquad \qquad \supset_{h} I(e_{Y_i}, \sigma_{Y_i}) \\[6pt]
(2) \rightarrow & \bar{Y}_i = \bar{A}_1 \bar{e}_{i1} + \bar{A}_2 \underline{e}_{i2} + \cdots + \bar{A}_J \underline{e}_{iJ} \\
& \qquad \qquad \qquad \qquad \supset_{h} I(e_{Y_i}, \sigma_{Y_i}) \\[6pt]
(3) \rightarrow & \bar{Y}_i = \bar{A}_1 \underline{e}_{i1} + \bar{A}_2 \bar{e}_{i2} + \cdots + \bar{A}_J \underline{e}_{iJ} \\
& \qquad \qquad \qquad \qquad \supset_{h} I(e_{Y_i}, \sigma_{Y_i}) \\[6pt]
\vdots & \qquad \qquad \qquad \vdots \\[6pt]
(2^J) \rightarrow & \bar{Y}_i = \bar{A}_1 \bar{e}_{i1} + \bar{A}_2 \bar{e}_{i2} + \cdots + \bar{A}_J \bar{e}_{iJ} \\
& \qquad \qquad \qquad \qquad \supset_{h} I(e_{Y_i}, \sigma_{Y_i})
\end{array}
\right\}
\end{equation}

It is an exhaustive way to solve the fuzzy random linear regression model (\ref{eq:2.10}). This problem consumes a substantial amount of computing time, particularly as $J$ increases. For instance, with $10^3$ features ($J$) and $10^4$ samples ($N$), the linear programming model (\ref{eq:2.10}) involves $2 \times 10^4 \times 2^{10^3}$ restrictions and $10^3$ non-negative restrictions. That is, $10^3 + 2 \times 10^4 \times 2^{10^3}$ (i.e. $J+2 \times 2^J \times N $)

which resulted in the linear programming having a huge number of constraints to solve. Therefore, we must employ a heuristic algorithm that was proposed in \cite{[A02]} to solve this model.
As explained in \cite{[A14]}, \cite{[A02]},
we can solve the problem by only $O(2^J)=2^{1,000}$ (\cite{[A02]}) in comparison to the complete solution above. Hence, the heuristic method is quite efficient in solving our T2F random linear regression model. 

\subsection{Heuristic Algorithm}\label{sec:3.5}

Note that, we give a new notations for $\bar{\textbf{A}}_j =[\underline{a}_j,\bar{a}_j]$, and step(n) denoted as $\bar{\textbf{A}}^{(n)}_j =[\underline{a}^{(n)}_j ,\bar{a}^{(n)}_j]$. 
In the Algorithm (\ref{algorithm6}), the multiplication 
formula $\bar{\textbf{A}}_j I(e_{X_{ij}}, \sigma_{X_{ij}})$, for $i  \in\{1,2,\dots,N\} $ and $j  \in\{1,2,\dots,J\} $, in equ. (\ref{(2.8)}) can be solved in various cases as shown in Table \ref{table2-2}.

At the level $h_0$, we denote an $\alpha$-level set of fuzzy degree for a structural attribute as follows: 
\begin{equation}\label{(2.11)}
(\mathbf{\bar{A}}_j)_{h_0} =[\underline{a}_j,\bar{a}_j] 
\end{equation}
for every $i$ and $j$, the signs of the confidence interval $I(e_{X_{ij}},\sigma_{X_{ij}})=[\underline{e}_{ij},\bar{e}_{ij}]$ and Eqn. (\ref{(2.11)}) restructure the interval of the product: 
\begin{align*}
(\mathbf{\bar{A}}_j \cdot I(e_{X_{ij}} ,\sigma_{X_{ij}} ))_{h_0} 
\end{align*}
is classified into several cases, as shown in Table \ref{table2-2}. 

As mentioned above, deriving analytical solutions for this problem is challenging, leading us to offer heuristic approaches. The procedure proposed is outlined in Algorithm (\ref{algorithm6}). 

\begin{table}[ht]
\caption{Various scenarios of the  proposed Product \cite{[A00]}}\label{table2-2}
\centering
\begin{tabular}{c|c|c}
\hline
Case & Condition & Result  \\ \cline{1-2}
Case I & $\bar{e}_{ij} \geq  \underline{e}_{ij} \geq 0$ & \\  \hline
I-a & $\bar{a}_j \geq  \underline{a}_j \geq 0$  & $(\bar{\textbf{A}}_j \cdot I(e_{X_{ij}}, \sigma_{X_{ij}}))_{h_0} = [\underline{a}_j \cdot \underline{e}_{ij}, \bar{a}_j \cdot \bar{e}_{ij}]$\\
I-b & $\bar{a}_j \geq 0 \geq \underline{a}_j$  & $(\bar{\textbf{A}}_j \cdot I(e_{X_{ij}} , \sigma_{X_{ij}}))_{h_0} = [\underline{a}_j \cdot \bar{e}_{ij}, \bar{a}_j \cdot \bar{e}_{ij}]$\\
I-c & $0 \geq \bar{a}_j \geq  \underline{a}_{ij}$ & $(\bar{\textbf{A}}_j \cdot I(e_{X_{ij}} , \sigma_{X_{ij}}))_{h_0} = [\underline{a}_j \cdot \bar{e}_{ij}, \bar{a}_j \cdot \underline{e}_{ij}]$\\
\hline
Case II & $0 \geq \bar{e}_{ij} \geq  \underline{e}_{ij}$ & \\  \hline
II-a & $\bar{a}_j\geq  \underline{a}_j \geq 0$  & $(\bar{\textbf{A}}_j \cdot I(e_{X_{ij}} , \sigma_{X_{ij}}))_{h_0} = [\bar{a}_j \cdot \underline{e}_{ij}, \underline{a}_j \cdot \bar{e}_{ij}]$\\
II-b & $\bar{a}_j \geq 0 \geq \underline{a}_j$  & $(\bar{\textbf{A}}_j \cdot I(e_{X_{ij}} , \sigma_{X_{ij}}))_{h_0} = [\bar{a}_j \cdot \underline{e}_{ij}, \underline{a}_j \cdot \underline{e}_{ij}]$\\
II-c & $0 \geq \bar{a}_j \geq  \underline{a}_j$  & $(\bar{\textbf{A}}_j \cdot I(e_{X_{ij}} , \sigma_{X_{ij}}))_{h_0} = [\bar{a}_j \cdot \bar{e}_{ij}, \underline{a}_j \cdot \underline{e}_{ij}]$\\
\hline
Case III & $\bar{e}_{ij} \geq  0 \geq \underline{e}_{ij}$ & \\  \hline
III-a & $\bar{a}_j \geq  \underline{a}_j \geq 0$  & $(\bar{\textbf{A}}_j \cdot I(e_{X_{ij}} , \sigma_{X_{ij}}))_{h_0} = [\bar{a}_j \cdot \underline{e}_{ij}, \bar{a}_j \cdot \bar{e}_{ij}]$\\
III-b & $0 \geq \bar{a}_j \geq  \underline{a}_j$  & $(\bar{\textbf{A}}_j \cdot I(e_{X_{ij}} , \sigma_{X_{ij}}))_{h_0} = [\underline{a}_j \cdot \bar{e}_{ij}, \underline{a}_j \cdot \underline{e}_{ij}]$\\
III-c & $\bar{a}_j \geq 0 \geq \underline{a}_j$ & $(\bar{\textbf{A}}_j \cdot I(e_{X_{ij}} , \sigma_{X_{ij}}))_{h_0} = [a^*_j \cdot e^*_{ij}, a^{**}_j \cdot e^{**}_{ij}]$\\
\hline 
\multicolumn{3}{p{8cm}}{*Note that 
$a^*_j \cdot e^*_{ij} = \min \{\underline{a}_j \cdot \bar{e}_{ij}\bar{a}_j \cdot \underline{e}_{ij}\}$, 
$a^{**}_j \cdot e^{**}_{ij} = \max \{\underline{a}_j \cdot \bar{e}_{ij}, \bar{a}_j \cdot \underline{e}_{ij} \}$} 
\end{tabular}
\end{table}

\begin{algorithm}
\caption{Heuristic Algorithm} \label{algorithm6}
\begin{algorithmic}
\State \textbf{Step 1 (Input):} Set the trial count $n \leftarrow 1$; define termination count $N_{\max*}$. Use $e_{ij}$ of attributes $A_j^{(n)}$, $j \in \{1,2,\dots,J\}$ for each sample $i \in \{1,2,\dots,N\}$, to find bounds $\underline{a}_j^{(n)}$ and $\bar{a}_j^{(n)}$ through solving the LP problem in equation (\ref{(2.8)}).

\State \textbf{Step 2 (Output):} Solve the LP to obtain initial bounds.

\State Determine $(\bar{A}_j \times I[e_{X_{ij}}, \sigma_{X_{ij}}])_{h^0}$ based on Cases I–III, using signs of $\bar{a}_j^{(n)}$. Compute $\underline{a}_j^{(n+1)}$, $\bar{a}_j^{(n+1)}$ and $A_j^{(n)}$ by solving LP to minimize fuzziness $J(\tilde{A})$, subject to constraints from (\ref{(2.8)}) and $\bar{A}_j \times I[e_{X_{ij}}, \sigma_{X_{ij}}]_{h_0}$.

\Procedure{Heuristics}{}
    \If{$\bar{a}_j^{(n+1)} \cdot \bar{a}_j^{(n)} \geq 0$ and $\underline{a}_j^{(n+1)} \cdot \underline{a}_j^{(n)} \geq 0$}
        \State Go to Step 5
    \Else
        \State $n \leftarrow n+1$
        \If{$n < N_{\max*}$}
            \State Go to Step 2
        \Else
            \State Go to Step 5
        \EndIf
    \EndIf
\EndProcedure

\State \textbf{Step 5 (Check Vertices):} For each vertex, check if it lies within the fuzzy linear regression lines. If inside or on either line, assign to set $S_1$; else assign to $S_2$.

\State \textbf{Step 6:} If $S_2$ is empty, proceed to Step 8. Otherwise, go to Step 7.

\State \textbf{Step 7:} Add points in $S_2$ as constraints in $LP^{(n)}$, resolve the LP.

\State \textbf{Step 8 (Terminate):} End the procedure.
\end{algorithmic}
\end{algorithm}

To address the computational challenges of exhaustive defuzzification, we shortly introduce how we sampling the linguistic terms into T2F data. In our scenario, we approximate the defuzzified value by randomly selecting a sample of embedded sets instead of considering all possible sets.
Continuous T2F sets have infinite embedded sets. The centroid values obtained through sampling are estimations. Discretization is a kind of sampling that selects specific points to represent the continuous set in a discrete manner.
Random selection of embedded sets is necessary due to their infinite nature. In grid-based discretization, secondary domain values lie at grid intersections. Creating an embedded set involves randomly selecting, for every primary domain value, a secondary domain value using a random function, ensuring equal probability. The selected embedded sets form a random sample, but the uniqueness among them is not guaranteed.

To determine the interval endpoints by using the Type-Resuction System (TRS) in the defuzzification process, we employ the Enhanced Iterative Algorithm (EIA) with Stop Condition, as outlined in \cite{[45]}. This iterative algorithm iteratively refines the interval boundaries until a specific stop condition is satisfied. Its purpose is to provide precise estimates of the TRS for interval endpoints.
In addition to the EIA, there are other methods available for defuzzify T2F into interval data. Two such methods are the Nie-Tan Method and the Greenfeld-Chiclana Collapsing Defuzzifier, as mentioned in reference \cite{[38]}. These methods operate by defuzzifying the $\alpha$-plane, which represents the primary membership grade, and then constructing a T1 TRS using the obtained defuzzified values. Typically, the defuzzified values generated by these methods are positioned approximately at the center of the interval.
Among the interval methods, the Collapsing Outward Right-Left (CORL) method, discussed in reference \cite{[42]}, has demonstrated high accuracy. It involves collapsing the interval outward from both the right and left endpoints. The CORL method aims to accurately represent the TRS interval by considering the distribution of the secondary membership grades.

In summary, these interval methods offer different approaches for defuzzification and determining the TRS interval endpoints. The selection of a method relies on considerations like the desired accuracy and computational complexity, and different methods may yield varying results depending on the specific scenario. Let us introduce our scenario in the following section.

\section{Real Scenario Cases and Numerical Example}\label{Sec:5} 
In this section, we consider a real scenario. A cosmetic company wants to plan a promotion with four categories of products for men and women: Basic Face Care (A), Face Cleaning (B), Cosmetics (C), Body Care (D). There are 3 products in each category. A fuzzy questionnaire was used to find out how satisfied men (10 people) and women (10 people) are that the products will reach the sales target. 
We also ask managers (10 people): Do you believe that each product will actually meet its new sales target based on your understanding and experiences? Participants can answer the questionnaires using linguistic terms such as "Always", "Frequently", "Often", "Sometimes", and "Seldom". The results of the survey are shown in Tables \ref{table3-1}, \ref{table3-2}, and \ref{table3}. The final goal of the company is to estimate the potential success (achieving sales targets) for each category of products, including men's and women's products.

\begin{table*}
\caption{View about the new sales target response by women}\label{table3-1}
\scriptsize 
\begin{tabular}{c|c|c|c|c|c|c|c|c|c|c|c}
\hline
\multirow{2}*{Productions} & Sales Volume & \multicolumn{10}{c}{Responses}\\  
\cline{3-12}
~ & (Millions) &  1 & 2 & 3 & 4 & 5 & 6 & 7 & 8 & 9 & 10 \\  
\hline
\multirow{3}*{A} & 5 & always & frequently & often & always & frequently & often & sometimes & seldom & always & often \\ 
~  & 5 & often &   frequently & often & sometimes&  seldom  & often & sometimes & seldom & frequently & often  \\
~  & 10 &  frequently &   often & frequently & frequently & often & sometimes & often &  frequently & seldom & often \\  \hline
\multirow{3}*{B} & 3 & seldom & often & sometimes &  often & frequently & often & sometimes & always & always & often  \\
~  & 4 & sometimes & often &  often & often & frequently & frequently & frequently & sometimes & seldom & often \\  
~ & 13 & always & often & always & often & frequently & frequently & frequently &  always & often & often  \\ \hline
\multirow{3}*{C} &  6 & frequently & often & frequently & frequently & often & sometimes & often & frequently & seldom & often \\  
~ & 6 & always & often & frequently & sometimes & often & frequently & seldom & always & often & often  \\
~ & 8 & often &  always & frequently & often & sometimes &  often & sometimes & often & frequently & often  \\ \hline
\multirow{3}*{D} & 7 & seldom & sometimes &  seldom & sometimes & frequently & seldom & always & seldom & frequently & often  \\
~  & 8 & always & always & frequently & always & sometimes & frequently & frequently & often & sometimes & often \\  
~  & 5 & sometimes &  seldom & sometimes &  often & frequently & often & often & sometimes & seldom & often   \\
\hline
\end{tabular}
\end{table*}

\begin{table*}
\caption{View about the new sales target response by men}\label{table3-2}
\scriptsize 
\begin{tabular}{c|c|c|c|c|c|c|c|c|c|c|c}
\hline
\multirow{2}*{Productions} & Sales Volume & \multicolumn{10}{c}{Responses}\\  
\cline{3-12}
~ & (Millions) &  1 & 2 & 3 & 4 & 5 & 6 & 7 & 8 & 9 & 10 \\  
\hline
\multirow{3}*{A} & 4 & often & frequently & often & sometimes & frequently & always & often & sometimes & often & always \\
~ &  4 & often & always & frequently & frequently & always & frequently & frequently & often & always & always \\ 
~ & 12 & always & frequently & frequently & often & always & often & sometimes & frequently & always & always \\ \hline
\multirow{3}*{B} & 5 & often & always & frequently & sometimes & seldom & frequently & often & sometimes & often & always \\
~ &   5 & always & often & sometimes & frequently & seldom & often & always & frequently & sometimes & always \\ 
~ & 10 & seldom & sometimes & seldom & often & frequently & sometimes & often & seldom & sometimes & always \\ \hline
\multirow{3}*{C} & 6 & sometimes & sometimes & sometimes & sometimes & seldom & sometimes & often & seldom & sometimes & always \\ 
~ &   7 & often & sometimes & frequently & sometimes & sometimes & always & often & sometimes & often & always \\
~ &   7 & always & always & always & frequently & frequently & frequently & often & frequently & frequently & always \\ \hline
\multirow{3}*{D} & 5 & often & often & seldom & seldom & frequently & often & sometimes & often & sometimes & always \\
~ &   6 & always & frequently & always & always & often & frequently & frequently & always & always & always \\ 
~ & 11 & frequently & often & sometimes & sometimes & sometimes & frequently & seldom & seldom & sometimes & always \\ 
 \hline
\end{tabular}
\end{table*}

\begin{table*}
\caption{View about the new sales target response by managers}\label{table3}
\scriptsize 
\begin{tabular}{c|c|c|c|c|c|c|c|c|c|c|c}
\hline
\multirow{2}*{Productions} & Sales Volume & \multicolumn{10}{c}{Responses}\\  
\cline{3-12}
~ & (Millions) &  1 & 2 & 3 & 4 & 5 & 6 & 7 & 8 & 9 & 10 \\  
\hline
\multirow{3}*{A} & 10	& Always	& Seldom	& Often	& Often	& Frequently	& Always	& Frequently	& Often	& Seldom	& Always\\
~ & 10	& Sometimes	& Frequently	& Sometimes	& Seldom & Seldom	& Always	& Frequently	& Seldom	& Frequently	& Frequently\\
~ & 20	& Seldom	& Seldom	& Always	& Often	& Seldom	& Seldom	& Always	& Always	& Seldom &	Sometimes\\
\hline										
\multirow{3}*{B} & 10	& Always	& Frequently & Often	& Always	& Frequently	& Sometimes	& Often	& Seldom	& Always	& Sometimes\\
~ & 14 & Often	& Sometimes	& Seldom	& Frequently	& Frequently	& Always	& Frequently	& Seldom	& Often	& Frequently\\
~ & 16 &	Frequently	& Sometimes	& Often	& Frequently & Often	& Frequently	& Sometimes	& Always	& Frequently	& Frequently\\
\hline										
\multirow{3}*{C} & 8	& Always	& Frequently	& Often	& Seldom	& Always	& Frequently	& Sometimes	& Often	& Frequently	& Seldom\\
~ & 12	& Frequently	& Always	& Often	& Sometimes &	Frequently	& Often	& Seldom	& Frequently	& Seldom	& Frequently\\
~ & 20 &	Often	& Frequently	& Sometimes	& Always &	Often	& Sometimes	& Seldom	& Always	& Often	& Seldom\\
\hline										
\multirow{3}*{D} & 12 &	Always	& Frequently	& Seldom	& Often &	Sometimes	& Always	& Frequently	& Often	& Seldom	& Frequently\\
~ & 12 &	Always	& Sometimes	& Frequently	& Often	& Seldom	& Always	& Often	& Frequently	& Sometimes &	Frequently\\
~ & 16	& Often	& Sometimes	& Always	& Sometimes	& Seldom	& Frequently	& Often	& Seldom	& Frequently &	Sometimes\\
\hline
\end{tabular}
\end{table*}

To capture these kinds of problems, we need to transfer the
linguistic terms into numerical values. Thus, we assigned numerical weights to each linquestic term as follows: (always, 0.9), (frequently, 0.7), (often, 0.5), (sometimes, 0.3), and (seldom, 0.1). This process results in a special case of the T2F set, as shown in Table \ref{table3-3}. Continuouslly, to capture the uncertainty in the secondary grade from Table \ref{table3-3}, the triangular fuzzy random data are considered. The Table \ref{table3-4} summarizes the input and output data for this scenario.

\begin{table}[ht]
\caption{Transfer the Linguistic Terms to T2F Set} \label{table3-3}
\centering
\begin{tabular}{cl}
\hline 
\hline
$X_{11}$= & (0.3/0.9+0.2/0.7+0.3/0.5+0.1/0.3+0.1/0.1)/5 \\
     & + (0.2/0.7+0.4/0.5+0.2/0.3+0.2/0.1)/5 \\
     & +(0.4/0.7 +0.4/0.5+0.1/0.3+0.1/0.1)/10 \\ \hline
$X_{21}$= & (0.2/0.9+0.1/0.7+0.4/0.5+0.2/0.3+0.1/0.1)/3 \\
     &+ (0.3/0.7+0.4/0.5+0.2/0.3+0.1/0.1)/4 \\
     & + (0.3/0.9+0.3/0.7+0.4/0.5)/13 \\ \hline
$X_{31}$= & (0.4/0.7+0.4/0.5+0.1/0.3+0.1/0.1)/6 \\
     & + (0.2/0.9+0.2/0.7+0.4/0.5+0.1/0.3+0.1/0.1)/6 \\
     & + (0.1/0.9+0.2/0.7+0.5/0.5+0.2/0.3)/8 \\ \hline
$X_{41}$= & (0.1/0.9+0.2/0.7+0.1/0.5+0.2/0.3+0.4/0.1)/7 \\
     & + (0.3/0.9+0.3/0.7+0.2/0.5+0.2/0.3)/8 \\
     & + (0.1/0.7+0.4/0.5+0.3/0.3+0.2/0.1)/5 \\ \hline
\hline
$X_{12}$= & (0.2/0.9+0.2/0.7+0.4/0.5+0.2/0.3)/4 \\
     & + (0.4/0.9+0.4/0.7+0.2/0.5)/4 \\
     & + (0.4/0.9+0.3/0.7+0.2/0.5+0.1/0.3)/12 \\ \hline     
$X_{22}$= & (0.2/0.9+0.2/0.7+0.3/0.5+0.2/0.3+0.1/0.1)/5 \\
     & + (0.3/0.9+0.2/0.7+0.2/0.5+0.2/0.3+0.1/0.1)/5  \\
    & + (0.1/0.9+0.1/0.7+0.2/0.5+0.3/0.3+0.3/0.1)/10 \\ \hline
$X_{32}$= & (0.1/0.9+0.1/0.5+0.6/0.3+0.2/0.1)/6 \\
     &+ (0.2/0.9+0.1/0.7+0.3/0.5+0.4/0.3)/7 \\ 
     & + (0.4/0.9+0.5/0.7+0.1/0.5)/7 \\ \hline
$X_{42}$= & (0.1/0.9+0.1/0.7+0.4/0.5+0.2/0.3+0.2/0.1)/5 \\
     &+ (0.6/0.9+0.3/0.7+0.1/0.5)/6 \\ 
    & + (0.1/0.9+0.2/0.7+0.1/0.5+0.4/0.3+0.2/0.1)/11 \\ \hline
\hline
$Y_1$= & (0.3/0.9+0.2/0.7+0.3/0.5+0.2/0.1)/10 \\
       & + (0.1/0.9+0.4/0.7+0.2/0.3+0.3/0.1)/10 \\
       & + (0.3/0.9+0.1/0.5+0.1/0.3+0.5/0.1)/20 \\  \hline
$Y_2$= & (0.3/0.9+0.2/0.7+0.2/0.5+0.2/0.3+0.1/0.1)/10 \\
       & + (0.1/0.9+0.4/0.7+0.2/0.5+0.1/0.3+0.2/0.1)/14 \\
       & + (0.1/0.9+0.5/0.7+0.2/0.5+0.2/0.3)/16  \\ \hline
$Y_3$= & (0.2/0.9+0.3/0.7+0.2/0.5+0.1/0.3+0.2/0.1)/8 \\ 
       & + (0.1/0.9+0.4/0.7+0.2/0.5+0.1/0.3+0.2/0.1)/12 \\
       & + (0.2/0.9+0.1/0.7+0.3/0.5+0.2/0.3+0.2/0.1)/20  \\ \hline
$Y_4$= &  (0.2/0.9+0.3/0.7+0.2/0.5+0.1/0.3+0.2/0.1)/12 \\
       & + (0.2/0.9+0.3/0.7+0.2/0.5+0.2/0.3+0.1/0.1)/12 \\
       & + (0.1/0.9+0.2/0.7+0.2/0.5+0.3/0.3+0.2/0.1)/16  \\ \hline
\hline
\end{tabular}
\end{table}

\begin{table}[ht]
\caption{Input Data and Output Data }\label{table3-4}
\centering
\begin{tabular}{cl}
\hline
1 & $X_{11}$ =(5,4,6)T ,0.39;(5,4,6)T ,0.27;(10,8,12)T ,0.34 \\
2 & $X_{21}$ =(3,2,4)T ,0.31;(4,3,5)T ,0.29;(13,11,15)T ,0.40 \\
3 & $X_{31}$ =(6,5,7)T ,0.32;(6,5,7)T ,0.35;(8,7,9)T ,0.33 \\
4 & $X_{41}$ =(7,6,8)T ,0.27;(8,7,9)T ,0.46;(5,7,6)T ,0.27 \\ \hline
1 & $X_{12}$ =(4,3,5)T ,0.29;(4,3,5)T ,0.37;(12,10,14)T ,0.35 \\
2 & $X_{22}$ =(5,4,6)T ,0.36;(5,4,6)T ,0.39;(10,8,12)T ,0.25 \\
3 & $X_{32}$ =(6,5,7)T ,0.21;(7,6,8)T ,0.32;(7,6,8)T ,0.47 \\
4 & $X_{42}$ =(5,4,6)T ,0.27;(6,5,7)T ,0.48;(11,9,13)T ,0.25 \\ \hline \hline
1 & $Y_1$ =(10,8,12)T ,0.40;(10,8,12)T ,0.32;(20,16,24)T ,0.28 \\
2 & $Y_2$ =(10,8,12)T ,0.34;(14,12,16)T ,0.31;(16,14,18)T ,0.35 \\
3 & $Y_3$ =(8,6,10)T ,0.35;(12,10,14)T ,0.34;(20,16,24)T ,0.31 \\
4 & $Y_4$ =(12,10,14)T ,0.35;(12,10,14)T ,0.36;(16,14,18)T ,0.29 \\ \hline
\end{tabular}
\end{table}

The results were used to plan a product mix suitable for both men ($X_{i1}$)and women ($X_{i2}$) in the marketing campaign. Our goal function, hence, is reduced to the following formula:
\begin{equation}\label{(goalF)}
\bar{Y}_i = \bar{A}_{1}I(e_{X_{i1}}, \sigma_{X_{i1}}) +\bar{A}_{2}I(e_{X_{i2}}, \sigma_{X_{i2}}), ~~~i= 1,2, 3,4
\end{equation}
where $I(e_{X_{ij}}, \sigma_{X_{ij}})$, for $j = 1, 2$, are the one-sigma confidence intervals shown in equation (\ref{(2.7)}). Since $N = 4$, $J = 2$, to solve this fuzzy random linear regression function, we consider a one-sigma confidence interval in this model.  We give the one-sigma confidence interval for input data and output data in Table \ref{table3-6}, and by using equation (\ref{(2.8)}), we have our T2F-LLR model as follows:
\begin{equation}\label{(3.0)}
\scriptstyle %\footnotesize % \tiny %
\left. 
\begin{array}{rl}  
\min & J(\bar{A}) =  \bar{a}_1^{(1)} -  \underline{a}_1^{(1)} +  \bar{a}_2^{(1)} -  \underline{a}_2^{(1)} \\
\text{subject to}\\
& \bar{a}_1^{(1)} \geq  \underline{a}_1^{(1)},  \bar{a}_2^{(1)} \geq  \underline{a}_2^{(1)} \\
& [\underline{a}_1^{(1)},\bar{a}_1^{(1)}] \cdot [4.19, 9.21] + [\underline{a}_2^{(1)},\bar{a}_2^{(1)}] \cdot [2.92, 10.76]  \\
& \supseteq [8.00,17.56]  \\  
& [\underline{a}_1^{(1)},\bar{a}_1^{(1)}] \cdot [2.53, 12.05]  + [\underline{a}_2^{(1)},\bar{a}_2^{(1)}] \cdot [3.95, 8.55]  \\
& \supseteq [10.54,16.14]  \\  
& [\underline{a}_1^{(1)},\bar{a}_1^{(1)}] \cdot [5.55, 7.77]  + [\underline{a}_2^{(1)},\bar{a}_2^{(1)}] \cdot [6.01, 7.51]  \\ 
& \supseteq [7.89,18.27]  \\  
& [\underline{a}_1^{(1)},\bar{a}_1^{(1)}] \cdot [5.52, 8.32]  + [\underline{a}_2^{(1)},\bar{a}_2^{(1)}] \cdot [4.50, 9.46]  \\ 
& \supseteq [10.98,15.34]  
\end{array} \right\}
\end{equation}

\begin{table*}
\caption{Input and Output Data in the Model}\label{table3-6}
\centering
%\scriptsize 
\begin{tabular}{c|cc|cc|cc}
\hline
i & $(e_{X_{i1}},\sigma_{X_{i1}})$ &  $I[e_{X_{i1}} ,\sigma_{X_{i1}} ]$ & $(e_{X_{i2}} ,\sigma_{X_{i2}} )$ & $I[(e_{X_{i2}} ,\sigma_{X_{i2}} ]$ & $(e_{Y_i} ,\sigma_{Y_i} )$ & $I[e_{Y_i} ,\sigma_{Y_{i}} ]$  \\ \hline
1 & (6.70,2.51) & [4.19,9.21] & (6.84,3.92) & [2.92,10.76] & (12.80,4.76) & [8.00,17.56]  \\
2 & (7.29,4.76) & [2.53,12.05] & (6.25,2.30) & [3.95,8.55] & (13.34,2.80) & [10.54,16.14] \\ 
3 & (6.66,1.11) & [5.55,7.77] & (6.79,0.78) & [6.01,7.57] & (13.08,5.19) & [7.89,18.27]  \\
4 & (6.92,1.40) & [5.52,8.32] & (6.98,2.48) & [4.50,9.46] & (13.16,2.18) & [10.98,15.34]  \\
\hline
\end{tabular}
\end{table*}

We consider the Heuristic Algorithm \ref{algorithm6} to solve the best parameter in the T2F-LLR model. 

We illustrate how to solve the model in equation (\ref{(3.0)}) as follows:

\textbf{Step 1:} Set up the initial confidence interval for
\begin{equation} \nonumber 
(\bar{A_j}^{(1)} \cdot I[e_{X_{ij}},\sigma_{X_{ij}}])_{h_0}=[\underline{a_j}^{(1)}\cdot e_{ij},~ \bar{a_j}^{(1)}\cdot e_{ij}],   
\end{equation}
where $i=1,2, 3, 4$, and $j=1, 2$.
We get each confidence interval $I[e_{X_{ij}},\sigma_{X_{ij}}]$ for $e_{ij}$ in Tabel \ref{table3-6}, and the model rearrange as follows:
\begin{equation}\label{(3.2)}
\scriptstyle %\footnotesize % \tiny %
\left. 
\begin{array}{rl}  
\min  & J(\bar{A})=  \bar{a}_1^{(1)} -  \underline{a}_1^{(1)} +    \bar{a}_2^{(1)} -  \underline{a}_2^{(1)} \\
\text{subject to}\\
& \bar{a}_1^{(1)} \geq  \underline{a}_1^{(1)},  \bar{a}_2^{(1)} \geq  \underline{a}_2^{(1)} \\
& [\underline{a}_1^{(1)}\cdot 6.70,\bar{a}_1^{(1)}\cdot 6.70] + [\underline{a}_2^{(1)}\cdot 6.84,\bar{a}_2^{(1)}\cdot 6.84]  \\
& \supseteq [8.00,17.56]  \\  
& [\underline{a}_1^{(1)}\cdot 7.29,\bar{a}_1^{(1)}\cdot 7.29] + [\underline{a}_2^{(1)}\cdot 6.25,\bar{a}_2^{(1)}\cdot 6.25] \\
& \supseteq [10.54,16.14]  \\  
& [\underline{a}_1^{(1)} \cdot 6.66,\bar{a}_1^{(1)}\cdot 6.66] + [\underline{a}_2^{(1)}\cdot 6.79,\bar{a}_2^{(1)}\cdot 6.79] \\ 
& \supseteq [7.89,18.27]  \\  
& [\underline{a}_1^{(1)} \cdot 6.92,\bar{a}_1^{(1)}\cdot 6.92] + [\underline{a}_2^{(1)}\cdot 6.98,\bar{a}_2^{(1)}\cdot 6.98] \\ 
& \supseteq [10.98,15.34]  
\end{array} \right\}
\end{equation}

It means that we need to solve the following linear programming model:
\begin{equation}\label{(38)}
\scriptstyle %\footnotesize % \tiny %
\left. 
\begin{array}{rl}  
\min  & J(\bar{A})=  \bar{a}_1^{(1)} -  \underline{a}_1^{(1)} +    \bar{a}_2^{(1)} -  \underline{a}_2^{(1)} \\
\text{subject to}\\
& \bar{a}_1^{(1)} \geq  \underline{a}_1^{(1)},  \bar{a}_2^{(1)} \geq  \underline{a}_2^{(1)} \\
& \underline{a}_1^{(1)}\cdot 6.70 + \underline{a}_2^{(1)}\cdot 6.84 \leq 8.00\\
& \bar{a}_1^{(1)}\cdot 6.70 + \bar{a}_2^{(1)}\cdot 6.84  \geq 17.56 \\
& \underline{a}_1^{(1)}\cdot 7.29 +\underline{a}_2^{(1)}\cdot 6.25 \leq 10.54\\
& \bar{a}_1^{(1)}\cdot 7.29 + \bar{a}_2^{(1)}\cdot 6.25 \geq 16.14 \\
& \underline{a}_1^{(1)} \cdot 6.66 + \underline{a}_2^{(1)}\cdot 6.79 \leq 7.89\\
& \bar{a}_1^{(1)}\cdot 6.66 + \bar{a}_2^{(1)}\cdot 6.79 \geq 18.27 \\  
& \underline{a}_1^{(1)} \cdot 6.92 + \underline{a}_2^{(1)}\cdot 6.98 \leq 10.98\\
& \bar{a}_1^{(1)}\cdot 6.92 + \bar{a}_2^{(1)}\cdot 6.98 \geq 15.34  
\end{array} \right\}
\end{equation}

We get the results that min $J(\bar{A})=1.53$, where $\bar{A_1}^{(1)}=[1.18, 1.18]$ and $\bar{A_2}^{(1)}=[0.00, 1.53]$. 

Now, we move to step 2.

\textbf{Step 2:} Since we got $\underline{a}_j^{(1)}$ and $\bar{a}_j^{(1)}$, for $j=1, 2$ are all nonnegative. Refering to Table \ref{table2-2}, we got new confidence interval from I-a result: 

\begin{equation} \nonumber 
(\bar{\textbf{A}}_j^{(2)} \cdot I[e_{X_{ij}}, \sigma_{X_{ij}}])_{h_0} = [\underline{a}_j^{(2)} \cdot \underline{e}_{ij}, \bar{a}_j^{(2)} \cdot \bar{e}_{ij}],
\end{equation}
where $i=1,2, 3, 4$, and $j=1, 2$. This change the model in equation (\ref{(3.0)}) to the next linear programming problems:

\begin{equation}\label{(39)}
\scriptstyle %\footnotesize
\left. \begin{array} {rl}
\min  & J(\bar{A})=  \bar{a}_1^{(2)} -  \underline{a}_1^{(2)} +    \bar{a}_2^{(2)} -  \underline{a}_2^{(2)} \\
\text{subject to } & \bar{a}_1^{(2)} \geq  \underline{a}_1^{(2)} \geq 0,  \bar{a}_2^{(2)} \geq  \underline{a}_2^{(2)} \geq 0\\
& [\underline{a}_1^{(2)} \cdot 4.19, \bar{a}_1^{(2)} \cdot 9.21]
+[\underline{a}_2^{(2)} \cdot 2.92, \bar{a}_2^{(2)} \cdot 10.76] \\\
& \supseteq [8.00,17.56] \\
& [\underline{a}_1^{(2)} \cdot 2.53, \bar{a}_1^{(2)} \cdot 12.05]
+[\underline{a}_2^{(2)} \cdot 3.95, \bar{a}_2^{(2)} \cdot 8.55] \\\
& \supseteq [10.54,16.14] \\
& [\underline{a}_1^{(2)} \cdot 5.55, \bar{a}_1^{(2)} \cdot 7.77]
+[\underline{a}_2^{(2)} \cdot 6.01, \bar{a}_2^{(2)} \cdot 7.57] \\\
& \supseteq  [7.89,18.27] \\
& [\underline{a}_1^{(2)} \cdot 5.52, \bar{a}_1^{(2)} \cdot 8.32]
+[\underline{a}_2^{(2)} \cdot 4.50, \bar{a}_2^{(2)} \cdot 9.46] \\\
& \supseteq  [10.98,15.34] \\ 
\end{array} \right\}
\end{equation}

That is, we need to solve the following linear programming model:
\begin{equation}\label{(40)}
\scriptstyle %\footnotesize % \tiny %
\left. 
\begin{array}{rl}  
\min  & J(\bar{A})=  \bar{a}_1^{(2)} -  \underline{a}_1^{(2)} +    \bar{a}_2^{(2)} -  \underline{a}_2^{(2)} \\
\text{subject to}\\
& \bar{a}_1^{(2)} \geq  \underline{a}_1^{(2)} \geq 0,  \bar{a}_2^{(2)} \geq  \underline{a}_2^{(2)} \geq 0 \\
& \underline{a}_1^{(2)}\cdot 4.19 + \underline{a}_2^{(2)}\cdot 2.92 \leq 8.00\\
& \bar{a}_1^{(2)}\cdot 9.21 + \bar{a}_2^{(2)}\cdot 10.76  \geq 17.56 \\
& \underline{a}_1^{(2)}\cdot 2.53 +\underline{a}_2^{(2)}\cdot 3.95 \leq 10.54\\
& \bar{a}_1^{(2)}\cdot 12.05 + \bar{a}_2^{(2)}\cdot 8.55 \geq 16.14 \\
& \underline{a}_1^{(2)} \cdot 5.55 + \underline{a}_2^{(2)}\cdot 6.01 \leq 7.89\\
& \bar{a}_1^{(2)}\cdot 7.77 + \bar{a}_2^{(2)}\cdot 7.57 \geq 18.27 \\  
& \underline{a}_1^{(2)} \cdot 5.52 + \underline{a}_2^{(2)}\cdot 4.50 \leq 10.98\\
& \bar{a}_1^{(2)}\cdot 8.32 + \bar{a}_2^{(2)}\cdot 9.46 \geq 15.34  
\end{array} \right\}
\end{equation}

We get the results that min $J(\bar{A})=0.93$, where $\bar{A_1}^{(2)}=[1.42, 2.35]$ and $\bar{A_2}^{(2)}=[0.00, 0.00]$. 

\textbf{Step 3:} Since we got $\underline{a}_j^{(2)} \cdot \underline{a}_j^{(1)} \geq 0$ and $\bar{a}_j^{(2)} \cdot \bar{a}_j^{(1)} \geq 0$, for $j=1, 2$, by Heuristic Algorithm \ref{algorithm6}, we move decision to {\bf Step 5}.

In this step, we need to check all vertices that are from equation (\ref{(39)}) and identified that which vertex doesn't include the range $I[e_{Y_i}, \sigma_{Y_i}]$. The examination is outlined as the equation (\ref{(41)}). Note that the $[\underline{a}_1^{(2)},\bar{a}_1^{(2)}]=[1.42, 2.35]$ and $[\underline{a}_2^{(2)},\bar{a}_2^{(2)}]=[0.00, 0.00]$.

\begin{equation}\label{(41)}
\begin{array}{rl}
(1) \rightarrow & [\underline{a}_1^{(2)},\bar{a}_1^{(2)}] \cdot 4.19 +[\underline{a}_2^{(2)},\bar{a}_2^{(2)}] \cdot 2.92 \supseteq [8.00,17.56] \\
(2) \rightarrow & [\underline{a}_1^{(2)},\bar{a}_1^{(2)}] \cdot 4.19 +[\underline{a}_2^{(2)},\bar{a}_2^{(2)}] \cdot  10.76   \supseteq [8.00,17.56]  \\
(3) \rightarrow & [\underline{a}_1^{(2)},\bar{a}_1^{(2)}] \cdot 9.21 +[\underline{a}_2^{(2)},\bar{a}_2^{(2)}] \cdot  2.92 \nsupseteq [8.00,17.56]  \\
(4) \rightarrow & [\underline{a}_1^{(2)},\bar{a}_1^{(2)}] \cdot 9.21 +[\underline{a}_2^{(2)},\bar{a}_2^{(2)}] \cdot  10.76 \nsupseteq [8.00,17.56] \\
(5) \rightarrow & [\underline{a}_1^{(2)},\bar{a}_1^{(2)}] \cdot 2.53 +[\underline{a}_2^{(2)},\bar{a}_2^{(2)}] \cdot  3.95 \supseteq [10.54,16.14] \\
(6) \rightarrow & [\underline{a}_1^{(2)},\bar{a}_1^{(2)}] \cdot 2.53 +[\underline{a}_2^{(2)},\bar{a}_2^{(2)}] \cdot  8.55 \supseteq [10.54,16.14]  \\
(7) \rightarrow & [\underline{a}_1^{(2)},\bar{a}_1^{(2)}] \cdot 12.05 +[\underline{a}_2^{(2)},\bar{a}_2^{(2)}] \cdot 3.95 \nsupseteq [10.54,16.14]  \\
(8) \rightarrow & [\underline{a}_1^{(2)},\bar{a}_1^{(2)}] \cdot 12.05 +[\underline{a}_2^{(2)},\bar{a}_2^{(2)}] \cdot  8.55 \nsupseteq [10.54,16.14]  \\
(9) \rightarrow & [\underline{a}_1^{(2)},\bar{a}_1^{(2)}] \cdot 5.55 +[\underline{a}_2^{(2)},\bar{a}_2^{(2)}] \cdot  6.01 \supseteq [7.89,18.27] \\
(10) \rightarrow & [\underline{a}_1^{(2)},\bar{a}_1^{(2)}] \cdot 5.55 +[\underline{a}_2^{(2)},\bar{a}_2^{(2)}] \cdot  7.57 \supseteq [7.89,18.27] \\
(11) \rightarrow & [\underline{a}_1^{(2)},\bar{a}_1^{(2)}] \cdot 7.77 +[\underline{a}_2^{(2)},\bar{a}_2^{(2)}] \cdot  6.01 \nsupseteq [7.89,18.27] \\
(12) \rightarrow & [\underline{a}_1^{(2)},\bar{a}_1^{(2)}] \cdot 7.77 +[\underline{a}_2^{(2)},\bar{a}_2^{(2)}] \cdot  7.57 \nsupseteq [7.89,18.27]\\
(13) \rightarrow & [\underline{a}_1^{(2)},\bar{a}_1^{(2)}] \cdot 5.52 +[\underline{a}_2^{(2)},\bar{a}_2^{(2)}] \cdot  4.50 \supseteq [10.98,15.34] \\
(14) \rightarrow & [\underline{a}_1^{(2)},\bar{a}_1^{(2)}] \cdot 5.52 +[\underline{a}_2^{(2)},\bar{a}_2^{(2)}] \cdot  9.46 \supseteq [10.98,15.34] \\
(15) \rightarrow & [\underline{a}_1^{(2)},\bar{a}_1^{(2)}] \cdot 8.32 +[\underline{a}_2^{(2)},\bar{a}_2^{(2)}] \cdot  4.50 \nsupseteq [10.98,15.34] \\
(16) \rightarrow & [\underline{a}_1^{(2)},\bar{a}_1^{(2)}] \cdot 8.32 +[\underline{a}_2^{(2)},\bar{a}_2^{(2)}] \cdot  9.46 \nsupseteq [10.98,15.34]
\end{array}
\end{equation}

Following this examination, we need to add 16 qualifications when the lower or upper ranges are not satisfied. That's what we need to add the set $S_2$ to the next linear programming.
\begin{equation}\nonumber 
S_2=\{(3), (4), (7), (8), (11), (12), (15), (16)\}
\end{equation}

Now, we move to {\bf Step 7}.

\textbf{Step 7:} Add all inequalities which in the Set $S_2$ into the model in equation (\ref{(40)}.
Hence, we get the next linear programming problem (\ref{(42)}):

\begin{equation}\label{(42)}
\scriptstyle %\footnotesize
\left. \begin{array}{rl}
\min & J(\bar{A}) = \bar{a}_1^{(2)} - \underline{a}_1^{(2)} + \bar{a}_1^{(2)} - \underline{a}_1^{(2)} \\
\text{subject to}\\
& \bar{a}_1^{(2)} \geq  \underline{a}_1^{(2)} \geq 0,  \bar{a}_2^{(2)} \geq  \underline{a}_2^{(2)} \geq 0 \\
& \underline{a}_1^{(2)}\cdot 4.19 + \underline{a}_2^{(2)}\cdot 2.92 \leq 8.00\\
& \bar{a}_1^{(2)}\cdot 9.21 + \bar{a}_2^{(2)}\cdot 10.76  \geq 17.56 \\
& \underline{a}_1^{(2)}\cdot 2.53 +\underline{a}_2^{(2)}\cdot 3.95 \leq 10.54\\
& \bar{a}_1^{(2)}\cdot 12.05 + \bar{a}_2^{(2)}\cdot 8.55 \geq 16.14 \\
& \underline{a}_1^{(2)} \cdot 5.55 + \underline{a}_2^{(2)}\cdot 6.01 \leq 7.89\\
& \bar{a}_1^{(2)}\cdot 7.77 + \bar{a}_2^{(2)}\cdot 7.57 \geq 18.27 \\  
& \underline{a}_1^{(2)} \cdot 5.52 + \underline{a}_2^{(2)}\cdot 4.50 \leq 10.98\\
& \bar{a}_1^{(2)}\cdot 8.32 + \bar{a}_2^{(2)}\cdot 9.46 \geq 15.34\\
(3) \rightarrow & \underline{a}_1^{(2)} \cdot 9.21 + \underline{a}_2^{(2)} \cdot 2.92 \leq 8.00 \\
& \bar{a}_1^{(2)} \cdot 9.21 + \bar{a}_2^{(2)} \cdot 2.92 \geq 17.56 \\
(4) \rightarrow & \underline{a}_1^{(2)} \cdot 9.21 + \underline{a}_2^{(2)} \cdot 10.76 \leq 8.00 \\
& \bar{a}_1^{(2)} \cdot 9.21 + \bar{a}_2^{(2)} \cdot 10.76 \geq 17.56 \\
(7) \rightarrow & \underline{a}_1^{(2)} \cdot 12.05 + \underline{a}_2^{(2)} \cdot 3.98 \leq 10.54 \\
& \bar{a}_1^{(2)} \cdot 12.05 + \bar{a}_2^{(2)} \cdot 3.98 \geq 16.14 \\
(8) \rightarrow & \underline{a}_1^{(2)} \cdot 12.05 + \underline{a}_2^{(2)} \cdot  8.55 \leq 10.54 \\
& \bar{a}_1^{(2)} \cdot 12.05 + \bar{a}_2^{(2)} \cdot 8.55 \geq 16.14 \\
(11) \rightarrow & \underline{a}_1^{(2)} \cdot 7.77 + \underline{a}_2^{(2)} \cdot 6.01 \leq 7.89 \\
& \bar{a}_1^{(2)} \cdot 7.77 + \bar{a}_2^{(2)} \cdot 6.01 \geq 18.27 \\
(12) \rightarrow & \underline{a}_1^{(2)} \cdot 7.77 + \underline{a}_2^{(2)} \cdot 7.57 \leq 7.89 \\
& \bar{a}_1^{(2)} \cdot 7.77 + \bar{a}_2^{(2)} \cdot 7.57 \geq 18.27 \\
(15) \rightarrow & \underline{a}_1^{(2)} \cdot 8.32 + \underline{a}_2^{(2)} \cdot 4.50 \leq 10.98 \\
& \bar{a}_1^{(2)} \cdot 8.32 + \bar{a}_2^{(2)} \cdot 4.50 \geq 15.34 \\
(16) \rightarrow & \underline{a}_1^{(2)} \cdot 8.32 + \underline{a}_2^{(2)} \cdot 9.46 \leq 10.98 \\
& \bar{a}_1^{(2)} \cdot 8.32 + \bar{a}_2^{(2)} \cdot 9.46 \geq 15.34
\end{array} \right\}
\end{equation}

We get the results that min $J(\bar{A})=0.93$, where $\bar{A_1}^{(2)}=[1.42, 2.35]$ and $\bar{A_2}^{(2)}=[0.00, 0.00]$. 

We can see that the results are the same as in step 2 and do not improve anymore. Hence, we terminate the algorithm and get our parameters for equation (\ref{(goalF)}) as follows:

\begin{equation}\nonumber
\begin{array}{l}
\bar{A}_1^l = 2.35, \bar{A}_1^r = 1.42, \bar{A}_2^l = \bar{A}_2^r = 0.00. 
\end{array}
\end{equation}

Hence, our T2F-LLR model with a confidence interval can be presented as follows: 

\begin{equation}\label{Result}
\begin{array}{rl}
\bar{Y}_i & = \bar{A}_1 I[e_{X_{i1}} , \sigma_{X_{i1}} ] + \bar{A}_2 I[e_{X_{i2}} , \sigma_{X_{i2}} ] \\ \\ 
          & = \displaystyle{(\frac{\bar{A}_l^1+\bar{A}_r^1}{2})}, \bar{A}_1^l, \bar{A}_1^r)I[e_{X_{i1}} ,\sigma_{X_{i1}} ]+
          \bar{A}_2 I[e_{X_{i2}} , \sigma_{X_{i2}} ] \\ \\
          &=[1.89,1.42,2.35]I[e_{X_{i1}} , \sigma_{X_{i1}} ] + 0.00 I[e_{X_{i2}} , \sigma_{X_{i2}} ]. 
\end{array}
\end{equation}

\begin{figure}[ht]
\centering
\includegraphics[width=4in, clip,keepaspectratio]{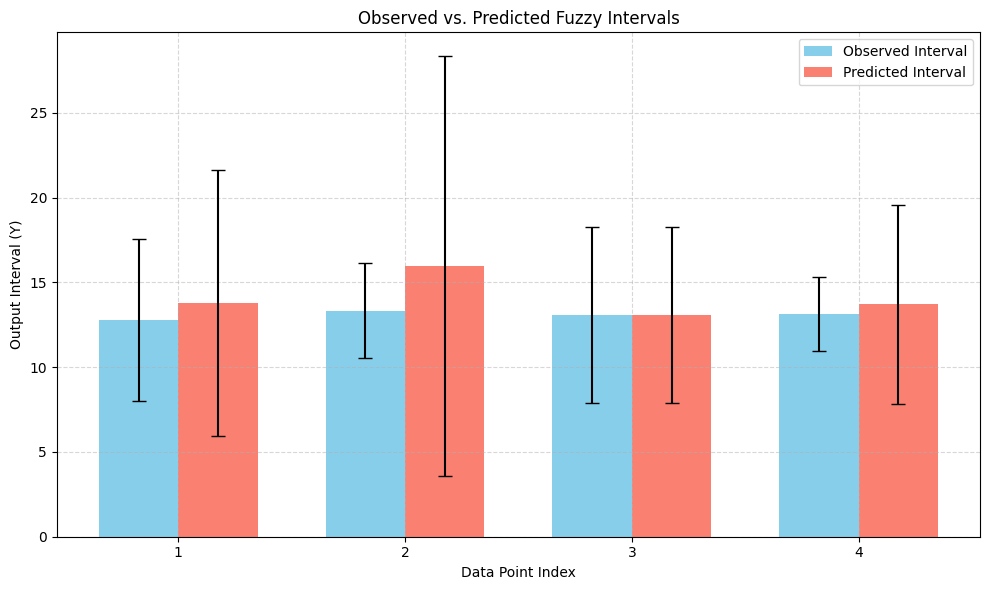}
\caption{Observed vs. Predicted Fuzzy Intervals}\label{OIPI}
\end{figure}

We give the observed confidence interval and the predicted confidence interval value in Table \ref{table3-5} and Figure \ref{OIPI} to check whether there are any differences in the four categories. We also provide the statistical accuracy metrics (Mean Absolute Percentage Error (MAPE) and Mean Squared Error (MSE) of weights in Table \ref{table3-5} by using the weight function, which was provided in the paper Lin {\it et. al.} in 2010 \cite{[lin2010]} that can deal with continuous fuzzy data for statistic test. Lin {\it et. al.} in 2016 \cite{[lin2016]} also present a statistical method which could distinguish the customer’s demand into different types, whereby fuzzy data is in consideration. Other influential work includes Lin {\it et. al.} in 2021 \cite{[lin2021]}. We give overall results in Table \ref{table3-5}.

\begin{figure}[ht]
\centering
\includegraphics[width=4in, clip,keepaspectratio]{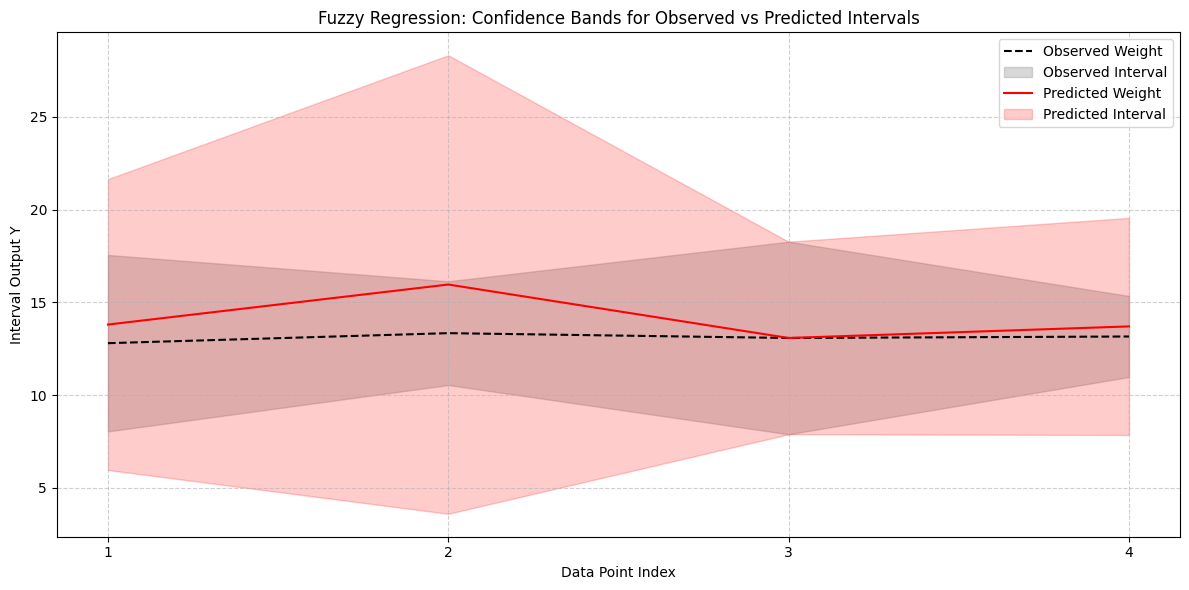}
\caption{Fuzzy Regression: Confidence Bands for Observed vs Predicted Intervals}\label{FRCIOIPI}
\end{figure}

From Table \ref{table3-5}, we can see that the results from the one-way ANOVA test (F-statistic = $2.6620$ and $p$-value = $0.1539$) show that there is no significant difference between Observed and Predicted Weights. If we want to see the significance of each production, the paired $t$-test (T-statistic = $1.8447$ and $p$-value = $0.1623$) shows that there is also no significant difference between Observed and Predicted Weights.
Overall, we can see the statistical accuracy metrics from MAPE, Product B (Face Cleaning) has a bigger MAPE (19.68\%) than other products; meanwhile, Product C (Cosmetics) almost has the same forecast result as the observed interval. It means that when a manager wants to use our model for making decisions about market strategy, they can focus more investment in cosmetics. Back to see our results in equation (\ref{Result}), we can see that the parameter $\bar{A_2}$ is almost equal $0.00$. It means that the model almost doesn't consider the input data from Table \ref{table3-2}, which is the questionnaire for men. For the above viewpoints, a manager can have a concept to make a market strategy to focus on "women" and "cosmetics".
For clear viewing, the optimal T2F-LinRLR model with confidence interval, the confidence bands for Observed vs Predicted Intervals are provided in Figure \ref{FRCIOIPI}.

\begin{table*}
\caption{Observed vs Predicted Fuzzy Intervals}\label{table3-5}
\centering
%\scriptsize 
%\tiny
%\small
\begin{tabular}{c|c|c|c|c|c}
\hline
Productions & i & Observed Interval & Predicted Interval & MAPE (Weight) & Difference (Weight)\\ \hline
A & 1 & [8.00, 17.56]  & [5.96, 21.66] & 8.06\% & 1.0035\\
B & 2 & [10.54, 16.14] & [3.60, 28.33] & 19.68\%& 2.6212\\
C & 3 & [7.89, 18.27]  & [7.89, 18.27] & 0\%    & 0.0000\\
D & 4 & [10.98, 15.34] & [7.85, 19.56] & 4.14\% & 0.5436\\
\hline
\multicolumn{2}{c|}{Statistical Accuracy Metrics} & \multicolumn{4}{l}{Coverage Rate: 100.00\% }\\
\multicolumn{2}{c|}{} & \multicolumn{4}{l}{MSE of Weight: 2.06 } \\
\multicolumn{2}{c|}{} & \multicolumn{4}{l}{MAPE of Weight: 7.97\% } \\
%\multicolumn{2}{c|}{} & \multicolumn{4}{l}{Mean Width Difference:8.16 }\\
\hline
\multicolumn{2}{c|}{One-Way ANOVA Test} & \multicolumn{4}{l}{F-statistic = 2.6620} \\
\multicolumn{2}{c|}{} & \multicolumn{4}{l}{$p$-value = 0.1539} \\
\multicolumn{2}{c|}{} & \multicolumn{4}{l}{Result: No significant difference between Observed and Predicted Weights.} \\
\hline
\multicolumn{2}{c|}{Paired $t$-test} & \multicolumn{4}{l}{T-statistic = 1.8447} \\
\multicolumn{2}{c|}{} & \multicolumn{4}{l}{$p$-value = 0.1623} \\
\multicolumn{2}{c|}{} & \multicolumn{4}{l}{Result: No significant difference between Observed and Predicted Weights overall.} \\
\hline
\end{tabular}
\end{table*}

\section{Conclusions and Future Works}\label{Sec:6} 
In this paper, we presented a novel heuristic algorithm for building a T2F-LLR model, which is designed to address complex decision-making problems involving both fuzziness and randomness. Motivated by a real-world scenario from the cosmetics industry, our approach effectively models linguistic assessments from both consumers and expert managers using T2F-LLR variables. We rigorously established a framework based on the one-sigma confidence interval under credibility theory to compute expected values and variances, providing a mathematically grounded method for interpreting fuzzy random input.
Our proposed heuristic algorithm offers a significant computational advantage by reducing the solving time without sacrificing accuracy, making it particularly suitable for non-meta datasets where rapid decision-making is critical. The experimental results demonstrated that our model achieves robust performance, with a low Mean Squared Error (MSE = 2.06) of weights, and all product-related variables were found to be statistically significant.

This study not only contributes to the theoretical foundation of fuzzy random regression modeling but also provides practical guidance for industry applications that involve linguistic and uncertain evaluations.
While this research establishes a solid basis for T2F-LLR modeling, several directions can be explored in future work, such as (1) Extension to Nonlinear or nonparametric regression frameworks, allowing for more flexibility in capturing complex relationships.
(2) The methodology can be adapted for use in domains such as healthcare, finance, and transportation, where linguistic assessments and uncertainty are prevalent.
(3) Further development of heuristic or metaheuristic optimization techniques (e.g., genetic algorithms, swarm intelligence) may enhance model scalability for larger datasets.

\bibliographystyle{unsrt}
\bibliography{T2FLLR0826}

\begin{IEEEbiography}[
{\includegraphics[width=1in,height=1.25in,clip,keepaspectratio]{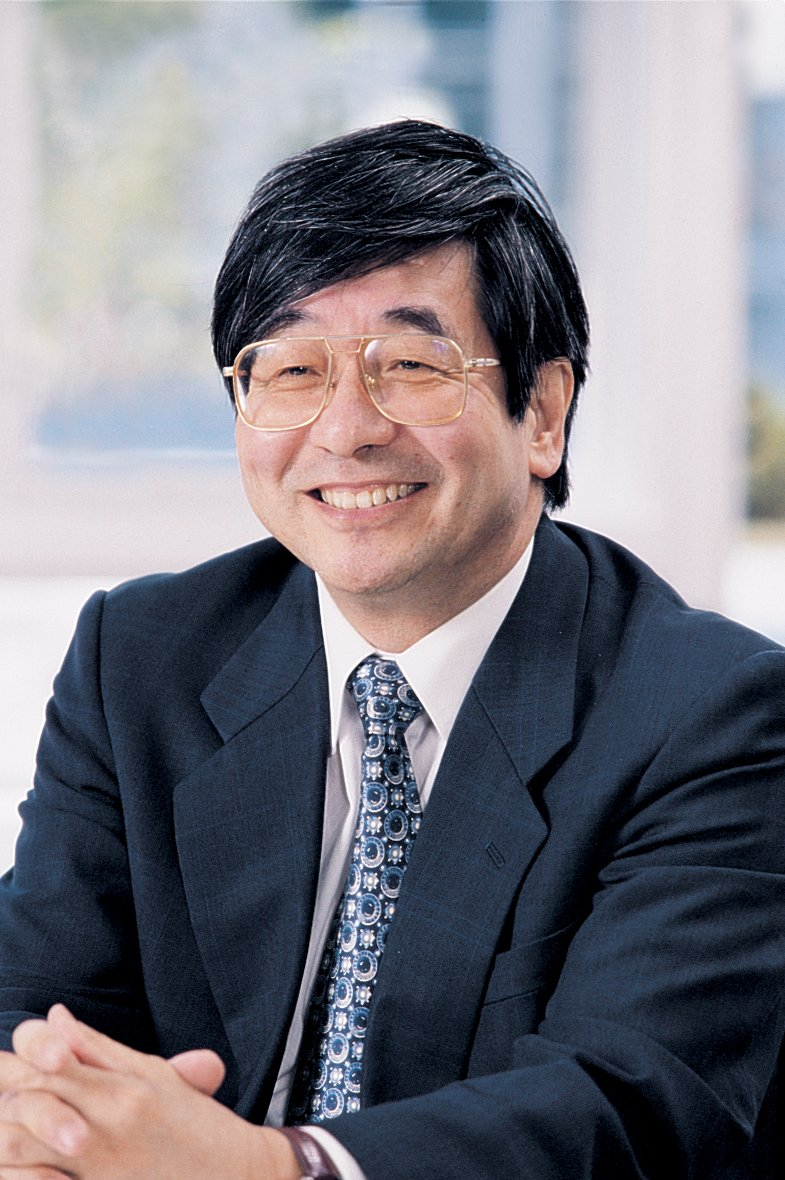}}]
{\bf Junzo Watada} (s' 80, M' 82, SM' 20, LSM' 24)
received the B.Sc. and M.Sc. degrees in electrical
engineering from Osaka City University, Osaka (now Osaka Metropolitan University),
Japan, and the Ph.D. degree from Osaka Prefecture University (now both are named Osaka Metropolitan University), Sakai, Japan. Until March 2016,
he was a Professor of management engineering,
knowledge engineering, and soft computing at the
Graduate School of Information, Production, and
Systems, Waseda University, Kitakyushu, Japan.
He is currently a specially appointed Professor, Fuculty of Data Science, Shimonoseki City University as well as  
since 2016, a Professor Emeritus. Waseda University. 
His research interests include big data analytics,
soft computing, tracking systems, knowledge engineering, and management
engineering. He is  a Life Senior Member, IEEE, a life fellow of the Japan Society for Fuzzy Theory and Intelligent Informatics as well as a life fellow of the Bio-Medical Fuzzy System Association.
He received the Henri Coanda Medal Award from Inventico, Romania,
in 2002; and the GH Asachi Medal from the Universitatea Technica GH
Asachi, IASI, Rumania, in 2006. 
\end{IEEEbiography}

\begin{IEEEbiography}[
{\includegraphics[width=1in,height=1.25in,clip,keepaspectratio]{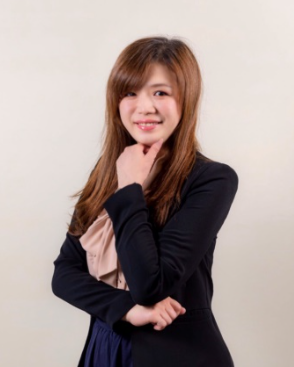}}]
{\bf Pei-Chun Lin} is currently an Associate Professor 
with the Department of Information
Engineering and Computer Science. She received her Ph.D. degree from the
Graduate School of Information, Production and Systems (IPS), Waseda University,
Japan. During her Ph.D. period, she built a series of statistical models for fuzzy data,
which was applied to decision-making systems in various fields. 
After she graduated
with her Ph.D., she worked as a researcher at IPS of Waseda University and
continued to specialize in the applications of fuzzy statistical models, and brought
fuzzy statistical models into artificial intelligence (AI) related applications. Dr. Lin not
only has good results in the research field but also serves as the editor and reviewer
of many top journals. At the same time, she also serves as the keynote speaker for
many international academic exchanges. 
Her research interests include Soft
Computing, Artificial Intelligence Computing, Robotics Computing, Statistical
Modeling, Cloud Computing, Big Data Analysis, etc.
\end{IEEEbiography}

\begin{IEEEbiography}[
{\includegraphics[width=1in,height=1.25in,clip,keepaspectratio]{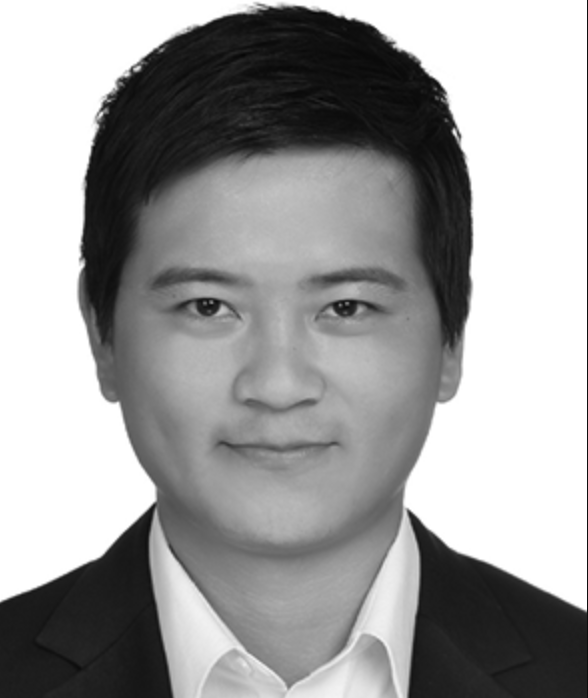}}]
{\bf Bo Wang}  (M’ 18) received the B.S. degree in software engineering from Southeast University, Nanjing, China, in 2007, and the M.S. and Ph.D. degrees
from the Graduate School of Information, Production
and Systems, Waseda University, Tokyo, Japan, in
2009 and 2012. 
He is currently an Associate Professor
with the School of Management and Engineering,
Nanjing University, Nanjing, China. He was a Research Assistant with Global COE Program, Waseda
University, Ministry of Education, Culture, Sports,
Science and Technology, Japan. His research interests
include power system planning, renewable generation forecasting, data-driven
decision-making, and artificial intelligence algorithms. 
He was a Special Research Fellow of the Japan Society for the Promotion of Science (JSPS). He
is a Committee Member of the Chinese Association of Automation Energy
Internet Committee, and a Committee Member of Systems Engineering Society
of JiangSu.
\end{IEEEbiography}

\begin{IEEEbiography}[
{\includegraphics[width=1in,height=1.25in,clip,keepaspectratio]{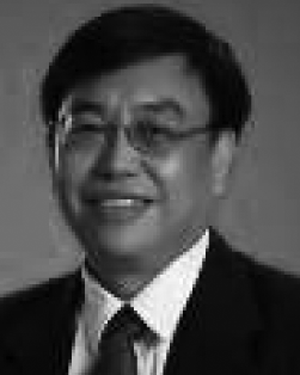}}]
{\bf Jeng-Shyang Pan} Jeng-Shyang Pan (Senior Member, IEEE) received the M.S. degree in communication engineering from National Chiao Tung University,
Hsinchu, Taiwan, in 1988, and the Ph.D. degree
in electrical engineering from the University of
Edinburgh, Edinburgh, U.K., in 1996.
He is currently a Professor with Nanjing University of Information Science \& Technology, Nanjing, China as well as 
the College of
Computer Science and Engineering, Shandong
University of Science and Technology, Qingdao,
China. His current research interests include
soft computing, information security, and signal
processing.
\end{IEEEbiography}

\begin{IEEEbiography}[
{\includegraphics[width=1in,height=1.25in,clip,keepaspectratio]{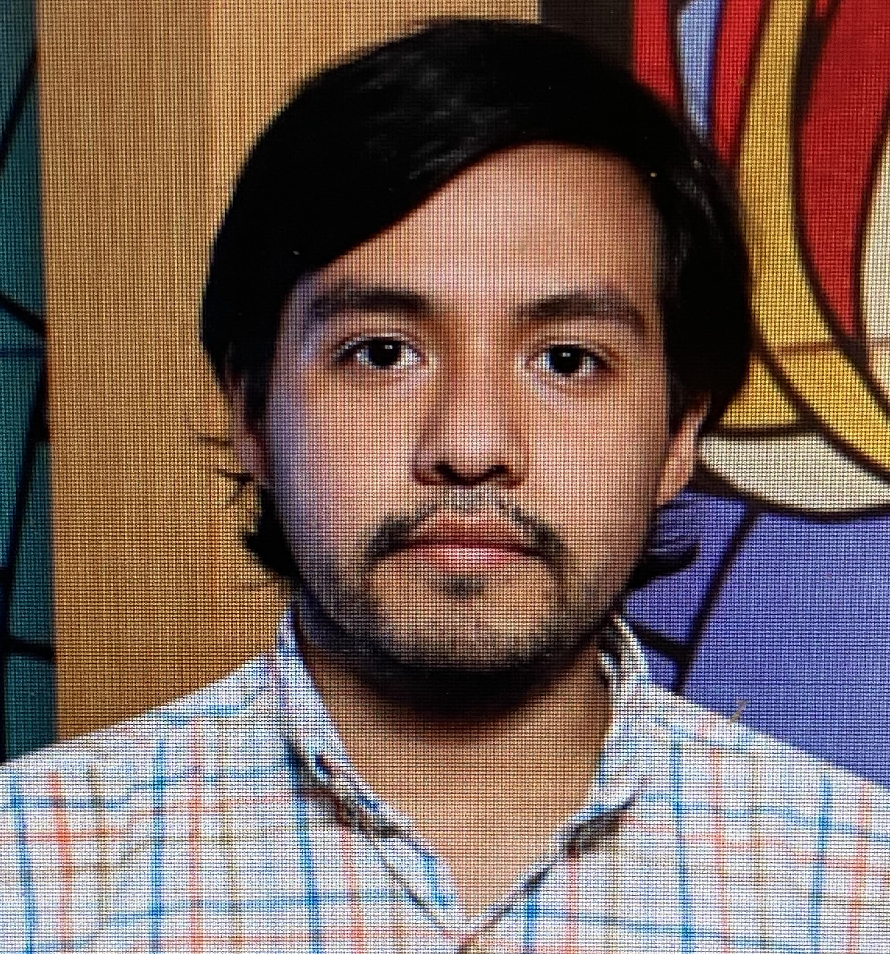}}]
{\bf Jos\'{e} Guadalupe Flores Mu\~{n}iz} received the Ph.D. Science Degree in Mathematics in 2020 from the Centro de Investigaci\'on de Ciencias F\'isico Matem\'aticas (CICFIM) of the Universidad Aut\'onoma de Nuevo Le\'on (UANL) in Mexico.
Since 2021, he has been working as a Professor and a Researcher with the CICFIM.
He is the author and the coauthor of one monograph, three book chapters, and 12 papers published in many prestigious journals, starting in 2016.
 His works are within the areas of game theory, bilevel programming, convex optimization, and conjectural variations.
\end{IEEEbiography}

\end{document}